\documentclass{article}
\usepackage{amsmath,amssymb,xspace}
\numberwithin{equation}{section} 

\advance\textwidth by 22pt\advance\oddsidemargin by-11pt
\advance\textheight by 37.5pt\advance\topmargin by-19pt

\newtheorem{theo}[equation]{Theorem}
\newtheorem{prop}[equation]{Proposition}
\newtheorem{cor}[equation]{Corollary}
\newtheorem{lem}[equation]{Lemma}

\def\beq{\begin{equation}}
\def\eeq{\end{equation}}
\def\p{\partial}
\def\G{\Gamma}
\def\g{\gamma}
\def\s{\sigma}
\def\z{\zeta}

\def\a{\alpha}
\def\b{\beta}
\def\e{\varepsilon}
\def\l{\lambda}
\def\f{\varphi}
\def\A{\mathcal A}
\def\V{\mathcal V}
\def\D{\mathcal D}
\def\F{\mathcal F}

\def\J{\mathcal J}
\def\L{\mathcal L}

\def\M{\mathcal M}

\def\O{\mathcal O}
\def\P{\mathcal P}

\def\wt{\widetilde}
\def\wh{\widehat}
\def\bC{\mathbb C}
\def\bZ{\mathbb Z}
\def\bN{\mathbb Z_{>0}}

\def\fm{\mathfrak m}

\def\Pic{\mathop{\rm Pic}\nolimits}
\def\Prym{\mathop{\rm Prym}\nolimits}
\def\dim{\mathop{\rm dim}\nolimits}
\def\res{\mathop{\rm res}\nolimits}

\def\BA{Baker-Akhiezer\xspace} 
\def\bH{\mathbb H}
\def\im{\mathop{\rm Im}\nolimits}
\def\Sp{\mathop{\rm Sp}\nolimits}

\begin{document}

\title{Soliton equations and the Riemann-Schottky problem}
\author{I. Krichever\thanks{Columbia University, New York, USA, and
Landau Institute for Theoretical Physics
and Kharkevich Institute for Problems of Information Transmission, Moscow, Russia, {\it Email address: \tt krichev@math.columbia.edu}\hfil\break
Research is supported in part by National Science
Foundation under the grant DMS-04-05519
and by The Ministry of Education and Science of the Russian Federation (contract 02.740.11.5194).}
\and T. Shiota\thanks{Kyoto University, Kyoto, Japan,
{\it Email address: \tt shiota@math.kyoto-u.ac.jp}}}
\date{}

\maketitle
\tableofcontents

\section{Introduction}
Novikov's conjecture on the Riemann-Schottky problem: {\it the Jacobians of smooth algebraic curves are precisely those indecomposable principally polarized abelian varieties (ppavs) whose theta-functions provide solutions to the Kadomtsev-Petviashvili (KP) equation}, was the first evidence of nowadays well-established fact: connections between the algebraic geometry and the modern theory of integrable systems is beneficial for both sides.

The purpose of this paper is twofold.
Our first goal is to present a proof of the strongest known characterization of a Jacobian variety in this direction:
{\it an indecomposable ppav $X$ is the Jacobian of a curve if and only if
its Kummer variety $K(X)$ has a trisecant line\/} \cite{kr-schot,kr-tri}.
We call this characterization {\it Welters' (trisecant) conjecture\/} after the work of
Welters \cite{wel1}. It was motivated by Novikov's conjecture and Gunning's celebrated theorem \cite{gun1}. The approach to its solution, proposed in \cite{kr-schot}, is general enough to be applicable to a variety of Riemann-Schottky-type problems. In \cite{kr-quad,kr-prym} it was used for a characterization of principally polarized Prym varieties. The latter problem is almost as old and famous as the Riemann-Schottky problem but is much harder. In some sense the Prym varieties may be geometrically the easiest-to-understand ppavs beyond Jacobians, and studying them may be a first step towards understanding the geometry of more general abelian varieties as well.

Our second and primary objective is to take this opportunity to elaborate on motivations underlining
the proposed solution of the Riemann-Schottky problem, to introduce a certain circle of ideas and methods, developed in the theory of soliton equations, and to convince the reader that they are algebro-geometric in nature, simple and universal enough to be included in the Handbook of moduli.
The results appeared in this article have already been published elsewhere.

\subsection*{Riemann-Schottky problem}

Let $\bH_g:=\{B\in M_g(\bC)\mid {}^tB=B,\ \im(B)>0\}$ be the Siegel upper half space.
For $B\in\bH_g$ let $\Lambda:=\Lambda_B:=\bZ^g+B\bZ^g$ and $X:=X_B:=\bC^g/\Lambda_B$.
Riemann's theta function
\beq\label{teta1}
\theta(z):=\theta(z,B):=\sum_{m\in\bZ^g}e^{2\pi i(m,z)+\pi i(m,Bm)},\quad
(m,z)=m_1z_1+\cdots+m_gz_g,
\eeq
is holomorphic and $\Lambda$-quasi\-periodic in $z\in\bC^g$,
so $\Theta:=\Theta_B:=\theta^{-1}(0)$ defines a divisor on $X$.
Moreover,
$(X,[\Theta])$ becomes a ppav, where
$[\Theta]$ denotes the algebraic equivalence class of $\Theta$.
Thus $\bH_g/\Sp(2g,\bZ)\simeq\A_g$, the moduli space of $g$-dimensional ppavs.
In what follows we may denote $(X,[\Theta])$ by $X$ for simplicity.
A ppav $(X,[\Theta])\in\A_g$ is said to be {\it indecomposable\/}
if $\Theta$ is irreducible, or equivalently\footnote{
since principal polarization means parallel translation is the only way
to deform $\Theta$, translating each component of $\Theta$ has
the same effect as translating $\Theta$ as a whole.
}
if there {\it do not\/} exist $(X_i,[\Theta_i])\in\A_{g_i}$ with $g_i>0$, $i=1$, 2,
such that $X=X_1\times X_2$ and $\Theta=\Theta_1\times X_2+X_1\times\Theta_2$.

Let $\M_g$ be the moduli space of nonsingular curves of genus $g$,
and let $J\colon\M_g\to\A_g$ be the Jacobi map, i.e., for $\G\in\M_g$,
$J(\G)$ is $\Pic^0(\G)$ with canonical polarization given by
$W_{g-1}=\{\L\in\Pic^{g-1}(\G)\mid h^0(\L)=h^1(\L)>0\}$
regarded as a divisor on $\Pic^0(\G)$, or more explicitly:
taking a symplectic basis $a_i$, $b_i$ ($i=1,\dots,g$) of
$H_1(\G,\bZ)$ and a basis $\omega_1$, \dots, $\omega_g$ of the space of
holomorphic 1-forms on $\G$ such that $\int_{a_i}\omega_j=\delta_{ij}$, we
define the {\it period matrix\/} and the {\it Jacobian variety\/} of $\G$ by
$$
B:=\left(\int_{b_i}\omega_j\right)\in\bH_g \quad\hbox{and}\quad
J(\G):=(X_B,[\Theta_B])\in\A_g\,,
$$
respectively. The latter is independent of the choice of $(a_i,b_i)$.

$J(\G)$ is indecomposable and the Jacobi map $J$ is injective (Torelli's theorem).
The {\it (Riemann-)Schottky problem\/} is the problem of characterizing the
Jacobi locus $\J_g:=J(\M_g)$ or its closure $\overline{\J_g}$ in $\A_g$.
For $g=2$,~$3$ the dimensions of $\M_g$ and $\A_g$ coincide,
and hence $\overline{\J_g}=\A_g$ by Torelli's theorem.
Since $\J_4$ is of codimension $1$ in $\A_4$, the case $g=4$ is the first nontrivial case of the Riemann-Schottky problem.

A nontrivial relation for the Thetanullwerte of a curve of genus~$4$
was obtained by F.~Schottky \cite{schottky} in 1888, giving a modular form which vanishes on
$\J_4$, and hence at least a {\it local\/} solution of the Riemann-Schottky
problem in $g=4$, i.e., $\overline{\J_4}$ is {\it an irreducible component
of\/} the zero locus $\mathcal S_4$ of the Schottky relation.
The irreducibility of $\mathcal S_4$ was proved by Igusa \cite{igusa}
in 1981, establishing $\overline{\J_4}=\mathcal S_4$, an effective answer
to the Riemann-Schottky problem in genus 4.

Generalization of the Schottky relation to a curve of higher genus, the so-called Schottky-Jung relations, formulated as a conjecture by Schottky and Jung \cite{schot-jung},
were proved by Farkas-Rauch \cite{far}. Later, van~Geemen \cite{geemen} proved that the Schottky-Jung relations give a local solution of the Riemann-Schottky problem. They do not give a global solution when $g>4$, since the variety they define has extra components already for $g=5$ (Donagi \cite{donagi2}).

More recent development on the Riemann-Schottky problem, as reviewed in
\cite{arb:expository,beauville,deb:expository}, includes a completely new
approach of Buser and Sarnak \cite{busersarnak} which provides an
effective way to characterize {\it non\/}-Jacobians.

\subsection*{Fay's trisecant formula and the KP equation}

Over more than 120 year-long history of the Riemann-Schottky problem,
quite a few geometric characterizations of the Jacobians have been obtained.
Following Mumford's review with a remark on Fay's trisecant formula
\cite{mum:c+j}, and the advent of soliton theory and Novikov's
conjecture \cite{kr1,kr2,mum}, much progress was made in the 1980s to
characterizing Jacobians and Pryms using Fay-like formulas and KP-like
equations. They are closely related to each other since Fay's formula,
written as a biliear equation for the Riemann theta function, follows
from a difference analogue of the {\it bilinear identity\/}\footnote{
Here $t=(t_1,t_2,\dots)$ and $t'=(t_1',t_2',\dots)$ are two sequences of
formal independent variables near zero, $k$ is a formal independent
variable near infinity,
$[k^{-1}]=(1/k,1/(2k^2),\penalty100 \dots,\penalty100 1/(nk^n),\dots)$,
and $\tau$, the so-called tau-function, is a
scalar-valued unknown function of the KP hierarchy. For a quasiperiodic
solution obtained from smooth curve $\G$ we have $\tau(t)=
e^{Q(t)}\theta(\sum t_iU_i+z,B(\G))$ for some quadratic form $Q(t)$,
vectors $U_i\in\bC^g$ and arbitrary $z\in\bC^g$.
Also, Fay's formula itself can in a sense be obtained from (\ref{KPBLI})
by specializing the time variables using the so-called Miwa variables.}
\beq\label{KPBLI}
\oint_{k=\infty}\tau(t-[k^{-1}])\tau(t'+[k^{-1}])e^{\sum(t_i-t'_i)k^i}dk=0
,
\eeq
which itself is equivalent to the KP hierarchy \cite{cherednik,sato}.
Equation (\ref{KPBLI}) can also be regarded as a generating function
for the Pl\"ucker relations for an infinite dimensional Grassmannian.

Compared with Igusa's work which studies the geometry of $\mathcal S_4$
and characterize the Jacobian {\it locus} $\J_4$, in this approach Fay-like
formulas or KP-like equations are used to (in a sense) construct the
curve $\G$ and thus characterize the Jacobian {\it varieties}.
Therefore this approach to the Riemann-Schottky problem is also
related to the Torelli theorem; however, the relation is only remote since
the conditions like Fay's formula and the KP equation contain extra
parameters like vector $U$ (and the lack of Prym-Torelli does not stop us
from studying the Prym-Schottky problem using the analogue of this
approach).



Let us first describe the trisecant formula in geometric terms.
The Kummer variety $K(X)$ of $X\in\A_g$ is the image of the Kummer map
\beq\label{kum}
K=K_X\colon X\ni z\longmapsto
\bigl(\Theta[\e,0](z)\mid \e\in((1/2)\bZ/\bZ)^g\bigr)\in \mathbb{CP}^{2^g-1}
\eeq
where
$\Theta[\e,0](z)=\theta[\e,0](2z,2B)$ are the level two theta-functions
with half-integer characteristics $\e\in((1/2)\bZ/\bZ)^g$, i.e., they equal
$\theta(2(z+B\e),2B)$ up to some exponential factor so that we have
\beq\label{ThetaQuad}
\theta(z+w)\theta(z-w)=\sum_{\e\in((1/2)\bZ/\bZ)^g}\Theta[\e,0](z)\Theta[\e,0](w)\,.
\eeq
We have $K(-z)=K(z)$ and $K(X)\simeq X/\{\pm1\}$.

A {\it trisecant\/} of the Kummer variety is a projective line which meets
$K(X)$ at three points.
{\it Fay's trisecant formula\/} states that if $X=J(\G)$, then $K(X)$ has a
family of trisecants parametrized by 4 points $A_i$, $1\le i\le4$, on $\G$.
Namely, identifying a point on $\G$ with its image under the Abel-Jacobi map
$\G\to\Pic^1(\G)$ and taking $r\in \Pic^{-1}(\G)$ such that
$2r=A_4-A_1-A_2-A_3$, we have:
\beq\label{fay:K}
K(r+A_1),\quad K(r+A_2)\quad\hbox{and}\quad K(r+A_3)\quad\hbox{are collinear,}
\eeq
i.e.,
\begin{gather*}
K\biggl(\frac{A_4+A_1-A_2-A_3}2\biggr),\quad
K\biggl(\frac{A_4-A_1+A_2-A_3}2\biggr)\quad\hbox{and}\\
K\biggl(\frac{A_4-A_1-A_2+A_3}2\biggr)\quad\hbox{are collinear}
\end{gather*}
if we take the three occurrences of ``division by 2'' consistent with each
other.
In what follows,
the same remark applies if division by 2 in $X$ appears more than once
in one formula, as in Theorems~\ref{th1.3}, \ref{theoremPrym}.

Since we have $K(-z)=K(z)$, condition (\ref{fay:K})
is symmetric in all the $A_i$'s. However, in its proof as well as its
applications the four points tend to play different roles.
E.g., fixing the 3 points $A_1$, $A_2$, $A_3$ we may regard it as
a one-parameter family of trisecants parametrized by $A_4$ or $r$.
Now drop the assumptions that $X=J(\G)$ and $A_i\in\G\subset X$:
suppose $X$ is a ppav such that (\ref{fay:K}) holds for some
$A_1$, $A_2$, $A_3\in X$ and infinitely many (hence a
one-parameter family of) $r\in X$.
Gunning proved in \cite{gun1} that, under certain non\-degeneracy
conditions, $X$ is then a Jacobian.


Gunning's work was extended by Welters who proved that a Jacobian variety can
be characterized by the existence of a formal one-parameter family of flexes of the
Kummer variety \cite{wel}. A flex of the Kummer variety
is a projective line which is tangent to $K(X)$ at some point up to order 2.
It is a limiting case of trisecants when the three intersection points come together.

In \cite{arb-decon} Arbarello and De~Concini showed that the assumption in
Welters' characterization is equivalent to a singly infinite sequence of
partial differential equations contained in the KP hierarchy, and proved
that only a first finite number of equations in the sequence are sufficient,
by giving an explicit bound for the number of
equations, $N=[(3/2)^gg!]$, based on the degree of $K(X)$.

\subsection*{Novikov's conjecture}

The second author's answer to Novikov's conjecture \cite{shiota} illustrated
how the soliton theory itself can provide natural, useful algebraic tools as
well as powerful analytic tools to study the Riemann-Schottky problem,
as immediately noticed by van der Geer \cite{vdG}, when only an early
version of \cite{shiota} was available:

An algebraic argument based on earlier results of Burchnall, Chaundy
and the first author \cite{ch,kr1,kr2} characterizes the Jacobians
using a commutative ring $R$ of ordinary differential operators associated
to a solution of the KP hierarchy.  A simple counting argument then shows
that only the first $2g+1$ time evolutions in the hierarchy are needed
to obtain $R$.
Indeed, suppose $X=\bC^g/\Lambda$ appears as an orbit of the first
$2g+1$ KP flows represented by a ``linear motion''
$\phi\colon\bC^{2g+1}\to
\bC^g$ followed by the projection $\bC^g\to X$.
Then $K:=\ker\phi$ is $(g+1)$-dimensional, and
if $(c_i)\in K$ then $\sum_ic_i\p\L/\p t_i=0$, hence
by the definition of the KP hierarchy $Q=\sum_ic_iP_i$ commutes with $\L$.
Any two such $Q$'s commute with each other \cite{schur}, so the $\bC$-algebra
$R'$ generated by all such $Q$'s is commutative.  A simple counting shows that
$R$ contains an ordinary differential operator of every order $n\ge2g+2$,
which implies that $R'$ is maximally commutative and hence $R'=R$,
from the way of constructing it.
Applying Burchnall et al's theory to $R$ to recover the spectral
curve $\G$ etc., we observe that $X\simeq J(\G)$.
The $2g+1$ KP flows yield a finite number of differential equations for
the Riemann theta function $\theta$ of $X$, to characterize a Jacobian.
As for the number of equations, an easy estimate shows that $4g^2$ is
enough, although more careful argument should yield a better bound.
Note that this is much smaller than Arbarello et al's estimate.


The analytic tools comes into play when one studies Novikov's conjecture,
that just the first equation ($N=1$!) of the hierarchy,
i.e., the KP equation (\ref{kp}), suffices to characterize the Jacobians:
in \cite{shiota} various tools obtained from analytic considerations on
the KP equation and family of its solutions were combined with the
algebraic arguments explained above to prove the conjecture.
Even Arbarello and De~Concini's geometric re-proof of
Novikov's conjecture \cite{arbarello} used the hardest analytic ingredient
of \cite{shiota} as it is, since it had no geometric alternative until
Marini's work \cite{mar} in 1998.
Analytic tools are also essential in the proofs of Welters' conjecture and
its Prym analogue presented in this paper, as condition (C) in each of
Theorems \ref{th1.1}, \ref{th1.2}, \ref{th1.3}, \ref{theoremPrym}.
Note that (\ref{cm50}), from which condition (C) in Theorem~\ref{th1.1}
follow, comes from a generalization of Calogero-Moser system.

Novikov's conjecture does not give an effective solution of the
Riemann-Schottky problem by itself: since it states that $X$ is a Jacobian
if and only if
$$
u=-2(\p_x^2\ln\theta(Ux+Vy+Wt+Z) + c)
$$
satisfies (\ref{kp}) for {\it some} $U$, $V$, $W$ and $c$,
we must eliminate those constants from (\ref{kp}) in order to
obtain an effective solution. It is hard to do this explicitly.

\subsection*{Welters' conjecture}

Novikov's conjecture is equivalent to the statement that the Jacobians are
characterized by the existence of length $3$ formal jet of flexes.
In \cite{wel1} Welters formulated the question: {\it if the Kummer variety $K(X)$ has {\it one\/}
trisecant, does it follow that $X$ is a Jacobian ?} In fact, there are
three particular cases of the Welters conjecture, corresponding to three possible configurations of
the intersection points $(a,b,c)$ of $K(X)$ and the trisecant:
\begin{itemize}
\item[(i)] all three points coincide $(a=b=c)$;
\item[(ii)] two of them coincide $(a=b\neq c)$;
\item[(iii)] all three intersection points are distinct
$(a\neq b\neq c\neq a)$.
\end{itemize}
Of course the first two cases can be regarded as degenerations of the general case~(iii).
However, when the presense of only one trisecant is assumed, all three cases are
independent and require separate treatment. The proof of case~(i) of
Welters' conjecture was obtained by the first author in \cite{kr-schot}:

\begin{theo}\label{th1.1} An indecomposable principally polarized abelian variety $(X,\theta)$
is the Jacobian variety of a smooth algebraic curve of genus g if and only if there exist $g$-dimensional vectors
$U\neq 0, V,A $ , and constants $p$ and $E$ such that one of the following three equivalent conditions are satisfied:

$(A)$ the equality
\beq
\left(\p_y-\p_x^2+u\right)\psi=0\,,\label{lax0}
\eeq
where
\beq\label{u0}
u=-2\p_x^2 \ln \theta (Ux+Vy+Z),\ \ \ \ \psi={\theta(A+Ux+Vy+Z)\over \theta(Ux+Vy+Z)}\, e^{p\,x+E\,y},
\eeq
holds, for an arbitrary vector $Z$;

\medskip
$(B)$ for all theta characteristics $\e\in ({1\over 2}\bZ/\bZ)^g$
$$ 
\left(\p_V-\p_U^2-2p\,\p_U+(E-p^2)\right)\, \Theta[\e,0](A/2)=0
$$ 
(here and below $\p_U$, $\p_V$ are the derivatives along the vectors $U$ and $V$, respectively).

\medskip
$(C)$ on the theta-divisor $\Theta=\{Z\in X\,\mid\, \theta(Z)=0\}$
\beq\label{cm70}
[(\p_V\theta)^2-(\p_U^2\theta)^2]\p_U^2\theta
+2[\p_U^2\theta\p_U^3\theta-\p_V\theta\p_U\p_V\theta]\p_U\theta+
[\p_V^2\theta-\p_U^4\theta](\p_U\theta)^2=0\ ({\rm mod}\, \theta)
\eeq
\end{theo}
The direct substitution of the expression (\ref{u0}) in equation (\ref{lax0}) and the use of the addition formula for the Riemann theta-functions shows the equivalence of conditions $(A)$ and $(B)$ in the theorem.
Condition $(B)$ 
means that the image of the point $A/2$ under the Kummer map is an inflection point (case~(i) of Welters' conjecture).

Condition $(C)$ is the relation that is {\it really used\/} in the proof of the theorem. Formally it is weaker than the other two conditions because its derivation does not use an explicit form (\ref{u0}) of the solution $\psi$ of equation (\ref{lax0}), but requires only an existence of a meromorphic solution: consider a holomorphic function $\tau(x,y)$ of a complex variable $x$ depending smoothly on a parameter $y$, and assume that in a neighborhood of a simple zero $\eta(y)$ of function $\tau$
(that is, $\tau(\eta(y),y)=0$ and $\p_x\tau(\eta(y),y)\ne0$) equation (\ref{lax0})
with potential $u=-2\p_x^2\ln \tau$ has a meromorphic solution $\psi$.
Then the equation
\beq\label{cm50}
\ddot \eta=2w,\
\eeq
holds, where the ``dots'' denote derivatives in $y$, and $w$ is
the third coefficient of the Laurent expansion of the function $u$ at the point $\eta$, i.e.,
$$
u(x,y)=\frac{2}{(x-\eta(y))^2}+v(y)+w(y)(x-\eta(y))+\cdots.
$$
Equations (\ref{cm50}) was first derived  in \cite{flex} where the assertion
of the theorem was proved under the assumption\footnote{under different
additional assumptions the corresponding statement was proved in the
earlier works \cite{kr3,mar}} that the closure of the group in $X$
generated by $A$ coincides with $X$. Expanding the function $\theta$
in a neighborhood of a point $z\in\Theta:=\{z\mid\theta(z)=0\}$ such that
$\p_U\theta(z)\ne0$, and noting that
the latter condition holds on a dense subset of $\Theta$ since $B$ is indecomposable,
it is easy to see that equation (\ref{cm50}) is equivalent to (\ref{cm70}).

Equation (\ref{lax0}) is one of the two auxiliary linear problems for the KP equation. Namely,
the compatibility condition of (\ref{lax0}) and the second auxiliary linear equation
\beq\label{lax01}
\left(\p_t-\p_x^3+\frac{3}{2}u\p_x+w\right)\psi=0
\eeq
is equivalent to the KP equation \cite{Dr,ZS1}:
\begin{equation}\label{kp}
\frac{3}{4} u_{yy}=\frac{\p}{\p x} \biggl(u_t-\frac{1}{4}u_{xxx}-
\frac{3}{2} u u_x \biggr)\,.
\end{equation}
For the first author, the motivation to consider not the whole KP equation but just one of its auxiliary linear problem was his earlier work \cite{kr3} on the elliptic Calogero-Moser (CM) system, where it was observed for the first time that equation (\ref{lax0}) is all what one needs to construct the elliptic solutions of the KP equation. Moreover, the construction of the Lax representation with a {\it spectral parameter\/} and the
corresponding spectral curves of the elliptic CM system proposed in \cite{kr3}
can be regarded as an effective solution of the inverse problem: how to reconstruct the algebraic curve from the matrix $B$ if its Kummer variety admits one flex with the vector $U$ (in the assumption of the Theorem) which {\it spans\/} an elliptic curve in the abelian variety $X$.
Briefly, that solution of the reconstruction problem can be presented as follows:

{\it If the vector $U$ spans an elliptic curve $E\subset X$, then the equation
\beq\label{algcm}
\theta(Ux+Vy+Z)=0
\eeq
for a generic $Z$ has $g$ simple roots $x_i(y)$ depending on $y$ (they are just intersection points of the shifted elliptic curve $E+Vy+Z\subset X$ with the theta-divisor $\Theta\subset X$). These roots define $g\times g$ matrix $L(y,z)$ with entries  given by
\begin{equation}\label{cmlax}
 L_{ii}(t,z)={1\over 2}\ \dot x_i,\ \ L_{ij}=\Phi(x_i-x_j,z), \ \ i\neq j,
\end{equation}
where
\begin{equation}\label{phi}
 \Phi(x,z):={\sigma(z-x)\over \sigma(z) \sigma(x)} e^{\zeta(z)x},
\end{equation}
with $\zeta$ and $\sigma$ the standard Weierstrass functions.

The spectral curve ${\G}_{cm}$ of the CM system is the normalization
at the point $k=\infty,z=0$ of the closure in ${\mathbb P}^1\times E$ of the affine curve given in $\mathbb C\times(E\setminus0)$ by the characteristic equation
\beq\label{r}
 R(k,z)=\det (kI+L(y,z))=0\,.
\eeq
Under the assumptions of the theorem, the CM curve $\G_{cm}$ does not depend on $y$ and is the solution of the inverse problem}.

Without an assumption on $U$ the proof of Theorem \ref{th1.1} is much more complex and less effective.
The ultimate goal is to construct, under the assumption that the condition $(C)$ is satisfied, a ring
of commuting ordinary differential operators, because, as shown in \cite{ch}, a pair of commuting differential operators $L_1,L_2$ satisfies an algebraic relation $R(L_1,L_2)=0$. This is the key moment, when
an algebraic curve emerges in the proof. It then remains only to show that the corresponding curve is the solution of the inverse problem.

The first step in the proof is to introduce in the problem a formal spectral parameter.
It is analogous to the introduction of the spectral parameter in the Lax matrix for the elliptic CM system.
This parameter $k$ appears in the notion of a {\it formal wave solution}
of equation (\ref{lax0}).

The wave solution of (\ref{lax0}) is a solution of the form
\beq\label{ps}
\psi(x,y,k)=e^{kx+(k^2+b)y}\biggl(1+\sum_{s=1}^{\infty}\xi_s(x,y)\,k^{-s}\biggr)\,.
\eeq
The aim is to show that under the assumptions of the theorem there exists
a unique, up to multiplication by a constant factor $c(k)$,  formal wave solution such that
\beq\label{xi22}
\xi_s={\tau_s(Ux+Vy+Z,y)\over\theta (Ux+Vy+Z)}.
\eeq
where $\tau_s(Z,y)$, is an entire function of $Z$.

As it was stressed above, strictly speaking the KP equation and the KP hierarchy are not present
in the assumptions of the theorem, but the analytical difficulties in the construction of the formal wave solutions of (\ref{lax0}) can be traced back to those in the second author's proof \cite{shiota} of Novikov's conjecture.

The main idea of proof in \cite{shiota} is to
show that if $\tau_0=e^{cx^2/2}\theta(Ux+Vy+Wt+Z)$ satisfies the KP equation in
Hirota's form\footnote{
We define $P(D_x,\dots)f\cdot f:=
P(\p_{x'},\dots)(f(x+x',\dots)f(x-x',\dots))|_{x'=\cdots=0}$ for a
polynomial or a power series $P$; a Hirota equation is an equation of
the form $P(D_x,\dots)f\cdot f=0$; see \cite{sato,shiota}.}
$$
(D_x^4+3D_y^2-4D_xD_t)\tau_0\cdot\tau_0=0,
$$
so that $u=-2\p_x^2\tau_0$ satisfies the KP equation (\ref{kp}), then
it can be extended to a $\tau$-function of the KP {\it hierarchy},
as a {\it global} holomorphic function of the infinite number of variables
$t=(t_i)=(t_1, t_2, t_3,\dots)$, with $t_1=x$, $t_2=y$, $t_3=t$. Local existence of $\tau$ directly follows from
the KP equation. The global existence of the $\tau$-function
is crucial. The rest is a corollary of the KP theory and the theory of commuting
ordinary differential operators developed by Burchnall-Chaundy \cite{ch} and
the first author \cite{kr1,kr2}.

The core of the problem is that there is a homological obstruction for the global existence
of $\tau$. It is controlled by the cohomology group $H^1(\mathbb C^g\setminus \Sigma, \V)$,
where {\it singular locus} $\Sigma$ is defined as $\p_U$-invariant subset of the theta-divisor
$\Theta$ and $\V$ is the sheaf of $\partial_U$-invariant meromorphic
functions on $\mathbb{C}^g\setminus \Sigma$ with poles along $\Theta$.
The hardest part of \cite{shiota}, as clarified in \cite{arbarello}, is
the proof that the locus $\Sigma$ is empty
\footnote{The first author is grateful to Enrico Arbarello
for an explanation of these deep ideas and a crucial role of the singular locus $\Sigma$,
which helped him to focus on the heart of the problem.}.

The coefficients $\xi_s$ of the wave function are defined recurrently by
the equation $2\partial_U\xi_{s+1}=\p_y \xi_s-\partial_U^2\xi_s+u\xi_s$. It turned out that
equation (\ref{cm70}) in the  condition $(C)$ of the theorem are necessary and sufficient for the local existence
of meromorphic solutions. The global existence of $\xi_s$ is controlled by the same cohomology group $H^1(\mathbb C^g\setminus \Sigma, \V)$ as above. Fortunately, in the framework of our approach
 there is no need to prove directly that the bad locus is empty.
The first step is to construct certain wave solutions outside the bad locus.
We call them $\l$-periodic wave solutions.
They are defined uniquely up to $\p_U$-invariant factor.
The next step is to show that for each $Z\notin \Sigma$ the $\l$-periodic wave solution
is a common eigenfunction of a commutative ring $\A^Z$ of ordinary difference operators.
The coefficients of these operators are independent of ambiguities in the construction of $\psi$.
For the generic $Z$ the ring $\A^Z$ is maximal
and the corresponding spectral curve $\G$ is $Z$-independent. The correspondence
$j\colon Z\longmapsto \A^Z$ and the results of the works \cite{ch,kr1,kr2,mum}, where
a theory of rank 1 commutative rings of differential operators was developed,
allows us to make the next crucial step and prove the global existence of the wave function.
Namely, on $(X\setminus\Sigma)$ the wave function can be globally defined as
the preimage $j^*\psi_{BA}$ under $j$ of the Baker-Akhiezer function on $\G$ and
then can be extended on $X$ by usual Hartogs' arguments. The global existence of the
wave function implies that $X$ contains an orbit of the
KP hierarchy, as an abelian subvariety. The orbit is isomorphic to the generalized
Jacobian $J(\G)={\rm Pic}^0(\G)$ of the spectral curve (\cite{shiota}).
Therefore, the generalized Jacobian
is compact. The compactness of $J(\G)$ implies that the spectral curve is smooth
and the correspondence $j$ extends by linearity and defines the isomorphism
$j\colon  X\to J(\G)$.

The proof of Welters' conjecture was completed in \cite{kr-tri}. First, here is
the theorem which treats case~(ii) of the conjecture:

\begin{theo}\label{th1.2} An indecomposable, principally polarized abelian variety $(X,\theta)$
is the Jacobian of a smooth curve of genus g if and only if
there exist non-zero $g$-dimensional vectors
$U\neq A \pmod \Lambda$, $V$, such that one of the following equivalent conditions
holds:

$(A)$  The differential-difference equation
\beq\label{laxd}
\left(\p_t-T+u(x,t)\right)\psi(x,t)=0, \ \ T=e^{\p_x}
\eeq
is satisfied for
\beq\label{u}
u=(T-1)v(x,t),\ \ v=-\p_t\ln\theta (xU+tV+Z)
\eeq
and
\beq\label{p}
\psi={\theta(A+xU+tV+Z)\over \theta(xU+tV+Z)}\, e^{xp+tE},
\eeq
where $p,E$ are constants and $Z$ is arbitrary.

\medskip
$(B)$ The equations
$$ 
\p_{V}\Theta[\e,0]\left((A-U)/2\right)-e^{p}\Theta[\e,0]\left((A+U)/2\right)
+E\Theta[\e,0]\left((A-U)/2\right)=0,
$$ 
are satisfied for all $\e\in({1\over 2}\bZ/\bZ)^g$. Here and below $\p_V$ is the constant
vector field on $\mathbb{C}^g$ corresponding to the vector $V$.

\medskip
$(C)$ The equation
\beq\label{cm7}
\p_V\left[\theta(Z+U)\,\theta(Z-U)\right]\p_V\theta(Z)=
\left[\theta(Z+U)\,\theta(Z-U)\right]\p^2_{VV}\theta(Z)\ ({\rm mod}\, \theta)
\eeq
is valid on the theta-divisor $\Theta=\{Z\in X\,\mid\, \theta(Z)=0\}$.
\end{theo}
Equation (\ref{laxd}) is one of the two auxiliary linear problems for the $2D$ Toda lattice
equation
\beq\label{2DT}
\p_\xi\p_\eta\varphi_n=e^{\varphi_{n-1}-\varphi_n}-e^{\varphi_{n}-\varphi_{n+1}},
\eeq
which can be regarded  as a partial discretization of the KP equation. The idea to use
it for the characterization of the Jacobians was motivated
by \cite{kr-schot} and the first author's earlier work with Zabrodin \cite{zab}, where
a connection of the theory of elliptic solutions of the $2D$ Toda lattice equations
and the theory of the elliptic Ruijsenaars-Schneider system was established.
In fact, Theorem~\ref{th1.2} in a slightly different form was proved in \cite{zab}
under the additional assumption that the vector $U$ {\it spans an elliptic curve} in $X$.

The equivalence of $(A)$ and $(B)$ is a direct corollary of the addition formula for
the theta-function. The statement $(B)$ is the second particular case of
the trisecant conjecture: the line in $\mathbb{CP}^{2^g-1}$ passing through the points
$K((A-U)/2)$ and $K((A+U)/2)$ of the Kummer variety is tangent
to $K(X)$ at the point $K((A-U)/2)$.

The affirmative answer to the third particular case, (iii), of Welters' conjecture is given by the following statement.

\begin{theo}\label{th1.3} An indecomposable, principally polarized abelian variety $(X,\theta)$
is the Jacobian of a smooth curve of genus g if and only if
there exist non-zero $g$-dimensional vectors
$U\neq V\neq A \neq U\, (\bmod \Lambda)$ such that one of the following equivalent
conditions holds:

$(A)$  The difference equation
\beq\label{laxdd}
\psi(m,n+1)=\psi(m+1,n)+u(m,n)\psi(m,n)
\eeq
is satisfied for
\beq\label{ud}
u(m,n)={\theta((m+1)U+(n+1)V+Z)\,\theta(mU+nV+Z)\over
\theta(mU+(n+1)V+Z)\,\theta((m+1)U+nV+Z)}
\eeq
and
\beq\label{pd}
\psi(m,n)={\theta(A+mU+nV+Z)\over \theta(mU+nV+Z)}\, e^{mp+nE},
\eeq
where $p,E$ are constants and $Z$ is arbitrary.

\medskip
$(B)$ The equations
$$ 
\Theta[\e,0]\left({A-U-V\over 2}\right)+e^{p}\Theta[\e,0]\left({A+U-V\over 2}\right)
=e^E\Theta[\e,0]\left({A+V-U\over 2}\right),
$$ 
are satisfied for all $\e\in({1\over 2}\bZ/\bZ)^g$.

\medskip
$(C)$ The equation
\beq\label{cm7d}
\theta(Z+U)\,\theta(Z-V)\,\theta(Z-U+V)+\theta(Z-U)\,\theta(Z+V)\,\theta(Z+U-V)=0 \pmod \theta
\eeq
is valid on the theta-divisor $\Theta=\{Z\in X\,\mid\, \theta(Z)=0\}$.
\end{theo}
Under the  assumption that the vector $U$ {\it spans an elliptic curve} in $X$,
Theorem~\ref{th1.3} was proved in \cite{bete}, where the connection of the
elliptic solutions of BDHE and, the so-called, elliptic nested Bethe Ansatz equations
was established.

Equation (\ref{laxdd}) is one of the two auxiliary linear problems for the so-called
bilinear discrete Hirota equation  (BDHE):
\beq
\tau _n (l+1,m)\tau _n (l,m+1)- \tau _n (l,m)\tau _n (l+1,m+1)+
\tau _{n+1} (l+1,m)\tau _{n-1} (l,m+1)=0\,
\label{BDHE}
\eeq
At the first glance all three nonlinear equation: the KP equation, the 2D Toda equation, and the BDHE equation,
look quite unlikely. But in the theory of integrable systems it is well-known that these
fundamental soliton equations are in intimate relation, similar to that between all three cases of the trisecant conjecture. Namely, the KP equation is as a continuous limit of the BDHE, and the 2D Toda equation
can be obtained in an intermediate step.

The structure of the statements of the last two theorems, and the structure of their proofs look almost literally
identical to that in Theorem~\ref{th1.1}. To some extend that is correct: in all cases the first step is to construct
the corresponding wave solution. The conditions $(C)$ in all three cases play the same role. They ensure the local
existence of the wave function. The key distinction between the differential and the difference cases arises at the next step. As it was mentioned above, in the case of differential equations a cohomological argument \cite[Lemma~12]{shiota}
can be applied to glue local solutions into a global one.
In the difference case there is no analog of the cohomological argument and we use
a different approach. Instead of {\it proving} the global existence
of solutions we, to some extend, {\it construct} them by defining first their residue on the theta-divisor.
It turns out that the residue is regular on $\Theta$ outside the {\it singular locus} $\Sigma$.
Surprisingly, it turns out that in the fully discrete case the proof of the statement that the singular locus is in fact empty can be obtained at much earlier stage than in the continuous or semi-continuous case. In part, it is due the drastic simplification in the fully discrete case of the corresponding equation on the theta-divisor
(compare (\ref{cm7d}) with (\ref{cm70})).

\subsection*{Structure of the article}
In the next section we introduce the basic concept of the algebro-geometric integration theory of soliton equation, that is the concept of the Baker-Akhiezer function,
which is defined by its analytic properties on an algebraic curve with fixed local coordinates at marked points.
The uniqueness of the Baker-Akhiezer function implies that it is a solution of certain linear differential equations. The existence of the Baker-Akhiezer function is proved by explicit theta-functional formula, which then
leads to explicit theta-functional formulae for the coefficients of the corresponding equations. That proves
``the only if'' part in all the theorems above.

In section 3, we introduce the KP hierarchy in Sato's form as a system of commuting flows on the space of formal pseudo\-differential operators. Because, the flows commute, the hierarchy can be reduced to the stationary points of one of the flows (or their linear combination). That is a reduction from a spatially two-dimensional system to a spatially one-dimensional system.\footnote{
Here the term ``spatially two-dimensional (resp.\ one-dimensional) system,''
also known as ``subsubholonomic (resp.\ subholonomic) system'' or
``$(2+1)$-d (resp.\ $(1+1)$-d) system,''
means the one whose ``general solution''
depends on functions of two variables (resp.\ one variable), or
equivalently, on doubly infinite (resp.\ singly infinite) sequences
of parameters in the formal power series set-up.
(The word ``space'' is associated to the notion of free parameters
because in an initial value problem of a partial differential equation
the free parameters for a solution are given by its initial data,
which are given on a ``space-like'' hypersurface.)
E.g., since initial data for the KP hierarchy, i.e.,
$\L|_{t=0}$ for $\L$ in (\ref{Lkp}), are given by a singly infinite sequence
of one-variable functions $\{v_s(x)\}_{s=1,2,\dots}$ or, by expanding
each $v_s(x)$ in a power series $v_s(x)=\sum_iv_{si}x^i$, a doubly infinite
sequence of parameters $\{v_{si}\}_{s=1,2,\dots;i=0,1,\dots}$,
the KP hierarchy is a ``spatially two-dimensional system.''
For $2\le n\in\bZ$ the $n$-reduction of the KP hierarchy (KdV if $n=2$,
Boussinesq if $n=3$, etc.) is defined by imposing the condition that
$\L^n$ is a differential operator.
Since, as an ordinary differential operator, $\L^n|_{t=0}$ depends on
finite number of one-variable functions and
hence on finite number of singly-infinite sequences, it is a ``spatially
one-dimensional system.''}
Under this reduction, the KP hierarchy defined first on ``a space'' of infinite number of functions of one variable (the coefficients of a pseudo\-differential operator) is equivalent to a system of commuting flows on the space of finite number of functions of one variable. For the case of stationary points of a linear combination of the first $n$ flows of the KP hierarchy these functions are coefficients of a differential operator $L_n$ of order $n$. One may take one step further and consider stationary
points of two commuting flows. It turns out that if the corresponding integers $n$ and $m$ are co-prime, then
the corresponding orbits of the whole hierarchy are finite-dimensional and can be identified with certain
subspaces of the finite-dimensional linear space of solutions to the system of ordinary differential equations:
\beq\label{comoper}
[L_n,L_m]=0,\ \ L_n=\p_x^n+\sum_{i=0}^{n-1} u_i(x)\p_x^{i},\ L_m=\p_x^n+\sum_{j=0}^{m-1} v_j(x)\p_x^{j}
\eeq
This is a setup explaining the role of commuting operators in the modern theory of integrable systems.

As a purely algebraic problem it was considered and partly solved in the remarkable
works of Burchnall and Chaundy \cite{ch} in the 1920s. They
proved that for any pair of such operators there exists  a polynomial
in two variables such that $R(L_n,L_m)=0$. Moreover, they proved that
if the orders $n$ and $m$ of these operators are co-prime, $(n,m)=1$,
and the algebraic curve $\G$ defined in $\mathbb C^2$ by equation $R(\lambda,\mu)=0$ is smooth, then
the commuting operators are uniquely defined by the curve and a set of $g$ points on $\G$, where $g$ is the genus of $\G$. In such a form, the solution of the problem is one of pure
classification: one set is  equivalent to the other. Even the attempt to obtain exact
formulae for the coefficients of commuting operators had not been made.
Baker proposed making the programme effective by looking at analytic properties
of the eigenfunction $\psi$. The Baker program was rejected by the authors of \cite{ch}
consciously (see the postscript of Baker's paper \cite{baker}) and all these results were forgotten for a long time.

The theory of commuting differential operators and its extension to the difference case is
presented in Section 4. The outline of the proof of the trisecant conjecture is in Section 5.
In Section 6 we present a solution of the characterization problem for Prym varieties
which was obtained by Grushevsky and the first author (\cite{kr-prym,kr-quad}). The last Section 7 is devoted to a
theory of abelian solutions of the soliton equation. The notion of such solutions was introduced by the authors in \cite{kr-shio,kr-shio1}, where it was shown that all of them are algebro-geometric.
The theory of abelian solutions can be regarded as an extension of the results above to the case of
non-principally polarized abelian varieties.

\section{The Baker-Akhiezer functions -- General scheme}

Let $\Gamma$ be a nonsingular algebraic curve of genus $g$ with $N$ marked points
$P_{\a}$ and fixed local parameters $k_{\a}^{-1}(Q)$ in  neighborhoods
of the marked points. The basic scalar {\it multi-point\/} and {\it multi-variable} Baker-Akhiezer function $\psi(t,Q)$ is a function of external parameters
\beq\label{times}
t=(t_{\a,i}),\ \a = 1,\ldots, N ; \ i=0,\ldots ;\ \  \sum_\a t_{\a,0}=0,
\eeq
only finite number of which is non-zero, and a point $Q\in \G$. For each set of the external parameters $t$ it is defined by its analytic properties on $\G$.

{\it Remark.} For the simplicity we will begin with the assumption that the variables $t_{\a,0}$ are integers, i.e.,
$t_{\a,0}\in \mathbb Z$.

\begin{lem} For any set of $g$ points  $\gamma_1,\ldots,\gamma_g$ in a
general position there exists a unique (up to constant factor
$c(t)$) function $\psi (t,Q)$,
such that:

(i) the function $\psi$ (as a function of the variable $Q\in \G$) is meromorphic everywhere except for the points $P_{\a}$ and
has at most simple poles at the points $\gamma_1,\ldots,\gamma_g$ ( if all
of them are distinct);

(ii) in a neighborhood of the point $P_{\a}$ the function $\psi$ has the
form
\beq
\psi (t,Q) =k_\a^{t_{\a,0}} \exp \biggl(\sum_{i=1}^{\infty} t_{\a ,i} k_{\a}^{i}
\biggr) \biggl( \sum_{s=0}^{\infty} \xi_{\a,s}(t) k_{\a}^{-s} \biggr),
\label{2.1}
\eeq
where $k_\a=k_\a (Q)$ is the reciprocal of a local parameter at $P_\a$, i.e.,
$k_\a^{-1}\in\frak m_{P_\a}\setminus m_{P_\a}^2$.
\end{lem}
 From the uniqueness of the Baker-Akhiezer function it follows that:

\begin{theo}\label{th2.1}
For each pair
$(\a,\,n>0)$ there exists a unique operator $L_{\a ,n}$ of the form
\beq
L _{\a ,n} = \p _{\a,1}^{n}
+ \sum_{j=0}^{n-1} u_{j}^{(\a ,n)}(t) \p_{\a ,1}^{j},
\label{2.2}
\eeq
(where $ \p _{\a,n} =\p / \p t _{\a ,n}$)
such that
\beq
\left(\p_{\a,n} - L_{\a,n}\right)\, \psi (t,Q) = 0 .
\label{2.3}
\eeq
\end{theo}
The idea of the proof of the theorems of this type proposed
in \cite{kr1}, \cite{kr2} is universal.

For any formal series of the form (\ref{2.1}) their exists a unique operator
$L_{\a ,n} $ of the form (\ref{2.2}) such that
\beq
\left(\p_{\a ,n} - L_{\a ,n} \right)\, \psi (t,Q) = O(k_\a^{-1})
\exp \,\biggl(\sum_{i=1}^{\infty} t_{\a ,i} k_{\a}^{i} \biggr) . \label{2.4}
\eeq
The coefficients of $L_{\a ,n} $ are universal differential polynomials with
respect to $\xi_{s,\a }$. They can be found after substitution of the
series (\ref{2.1}) into (\ref{2.4}).

It turns out that if the series (\ref{2.1}) is not formal but is an
expansion of the Baker-Akhiezer function in the neighborhood of $P_{\a}$ the
congruence (\ref{2.4}) becomes an equality. Indeed, let us consider the
function $\psi_{1}$
\beq
\psi_{1} = (\p_{\a ,n} - L_{\a ,n}) \psi (t,Q).
\label{2.5}
\eeq
It has the same analytic properties as $\psi$ except for the only one.
The expansion of this function in the neighborhood of $P_{\a}$ starts
from $O(k_\a^{-1})$.
 From the uniqueness of the Baker-Akhiezer function it follows that
$\psi_1 = 0 $ and the equality (\ref{2.3}) is proved.

\begin{cor} The operators $ L_{\a ,n}$ satisfy the compatibility
conditions
\beq
\bigl[ \p_{\a ,n} - L_{\a ,n} ,
\p_{\a ,m} - L_{\a ,m} \bigr] = 0 .\label{2.6}
\eeq
\end{cor}
{\bf Remark.} The equations (\ref{2.6}) are gauge invariant. For any function
$c(t)$ operators
\beq
\wt L_{\a ,n} = c L_{\a ,n} c^{-1} +
( \p_{\a ,n}c) c^{-1} \label{2.7}
\eeq
have the same form (\ref{2.2}) and satisfy the same operator equations
(\ref{2.6}). The gauge transformation (\ref{2.7}) corresponds to the gauge
transformation of the Baker-Akhiezer function
\beq
\wt\psi (t,Q) = c(t) \psi (t,Q) \label{2.7a}
\eeq
In addition to differential equations (\ref{2.3}) the \BA function satisfies an infinite system of
differential-difference equations. Recall that the discrete variables $t_{\a,0}$ are subject to the constraint
$\sum_\a t_{\a,0}=0$. Therefore, only the first $(N-1)$ of them are independent and $t_{N,0}=-\sum_{\a=1}^{N-1}t_{\a,0}$. Let us denote by $T_{\a},\ \ \a=1,\ldots,N-1,$ the operator that shifts
the arguments $t_{\a,0}\to t_{\a,0}+1$ and $t_{N,0}\to t_{N,0}-1$, respectively. For the sake of brevity in the formulation of the next theorem we introduce the operator $T_N=T_1^{-1}$.
\begin{theo}
For each pair
$(\a,\,n>0)$ there exists a unique operator $\wh L_{\a ,n}$ of the form
\beq
\wh L_{\a ,n} = T_{\a}^{n}+ \sum_{j=0}^{n-1} v_{j}^{(\a ,n)}(t) \,T_{\a}^{j},\ \
v_{0}^{(N ,n)}(t)=0.
\label{2.2d}
\eeq
such that
\beq
\left(\p_{\a,n} - \wh L_{\a,n}\right)\, \psi (t,Q) = 0 .
\label{2.3d}
\eeq
\end{theo}
The proof if identical to that in the differential case. The operators $\wh L_{\a,n}$ are defined by
congruence insuring that the resulting function satisfies all the condition of the \BA function plus vanishing
of one of the leading coefficients. After that the uniqueness of the \BA function implies that the congruence is in fact the equality.
\begin{cor} The operators $\wh L_{\a ,n}$ satisfy the compatibility
conditions
\beq
\bigl[ \p_{\a ,n} - \wh L_{\a ,n} ,
\p_{\a ,m} - \wh L_{\a ,m} \bigr] = 0 .\label{2.6d}
\eeq
\end{cor}
It should be emphasized that the algebro-geometric construction
is not a sort of abstract ``existence'' and ``uniqueness'' theorems. It
provides the explicit formulae for solutions in terms of the Riemann
theta-functions. They are the corollary of the explicit formula for the
Baker-Akhiezer function:
\begin{theo} The Baker-Akhiezer function is given by the formula
\beq
\psi(t,P)=c(t)\exp\left(\sum t_{\a,i}\Omega_{\a,i}(P)\right)
{\theta(A(P)+\sum U_{\a,i}t_{\a,i}+Z)
\over \theta(A(P)+Z)} , \label{2.101}
\eeq
Here the sum is taken over all the indices $(\a, i>0)$ and over the indices $(\a,0)$ with $\a=1,\ldots,N-1$, and:

a) $\Omega_{\a,i}(P)$ is the abelian integral,
$
\Omega_{\a,i}(P)=\int^P d\Omega_{\a,i},
$
corresponding to the unique normalized,
$
\oint_{a_k} d\Omega_{\a,i}=0,
$
meromorphic differential on
$\Gamma$, which for $i>0$ has the only pole of the form
$
d\Omega_{\a,i}=d\left(k_{\a}^i+O(1)\right)
$
at the marked point $P_{\a}$ and for $i=0$ has simple poles at the marked point $P_{\a}$ and $P_{N}$ with residues
$\pm 1$, respectively;

b) $2\pi iU_{\a,j}$ is the vector of $b$-periods of the differential
$d\Omega_{\a, j}$, i.e.,
$$
U_{\a,j}^k={1\over 2\pi i} \oint_{b_k} d\Omega_{\a,j};
$$

c) $A(P)$ is the Abel transform, i.e., a vector with the coordinates
$
A(P)=\int^P d\omega_k
$

d) $Z$ is an arbitrary vector (it corresponds to the divisor of poles of
Baker-Akhiezer function).
\end{theo}
Notice, that from the bilinear Riemann relations it follows that the expansion of the Abel transform near the marked point has the form
\beq\label{abelexpans}
A(P)=A(P_{\a})-\sum_{i=1}^\infty \frac{1}{i}U_{\a,i}k_{\a}^{-i}
\eeq

\subsubsection*{Example 1. One-point \BA function. KP hierarchy}

In the one-point case the Baker-Akhiezer function has an exponential
singularity at a single point $P_1$ and depends on a single set of variables $t_i=t_{1,i}$. Note that in this case there is no discrete variable, $t_{1,0}\equiv 0$.
Let us  choose the normalization of the Baker-Akhiezer function with the
help of the condition $\xi_{1,0}= 1 $, i.e., an expansion of $\psi$ in the neighborhood
of $P_1$ equals
\beq
\psi (t_1,t_2,\ldots,Q) =
\exp \biggl(\sum_{i=1}^{\infty} t_{i} k^{i} \biggr)
\biggl(1+ \sum_{s=1}^{\infty} \xi_{s}(t) k^{-s}\biggr).\label{2.8}
\eeq
Under this normalization (gauge) the corresponding operator $L_n$ has the form
\beq
L_n =\p_1^n +\sum_{i=0}^{n-2} u_i^{(n)} \p_1^i . \label{2.9}
\eeq
For example, for $n=2,3$ after redefinition $x=t_1$ we have
$L_2=\p_x^2-u, L_3=\p_x^3-\frac32u\p_x-w$
with
\beq
u(x,t_2,\ldots)=2\p_x \xi_1(x,t_2,\ldots), \label{2.10}
\eeq
Therefore, if we define  $y=t_2, t=t_3$, then
$ u(x,y,t,t_4,\ldots )$ satisfies the KP equation (\ref{kp}).

The normalization of the leading coefficient in (\ref{2.8}) defines the the function $c(t)$ in (\ref{2.101}).
That gives the following formula for the normalized one-point \BA function:
\beq
\psi(t,Q)=\exp\left(\sum t_{i}\Omega_{i}(P)\right)
{\theta(A(P)+\sum U_{i}t_{i}+Z) \,\theta(Z)
\over \theta(\sum U_{i}t_{i}+Z)\,\theta(A(P)+Z)} , \label{BA1}
\eeq
(shifting $Z$ if needed we may assumed that $A(P_1)=0$).
In order to get the explicit theta-functional form of the solution of the KP
equation it is enough to take the derivative of the first coefficient of the
expansion at the marked point of the ratio of theta-functions in the formula (\ref{BA1}).

Using (\ref{abelexpans}) we get the final formula for the algebro-geometric solutions of the KP hierarchy \cite{kr2}
\beq
u(t_1,t_2, \ldots ) = - 2 \p_1^2 \ln \theta (\sum_{i=1}^{\infty} U_i t_i +Z) +
\hbox{const}. \label{2.11}
\eeq

\subsubsection*{Example 2. Two-point \BA function. $2D$ Toda hierarchy}

In the two-point case the Baker-Akhiezer function has exponential
singularities at two points $P_\a, \a=1,2,$ and depends on two sets of continuous variables $t_{\a,i>0}$. In
addition it depends on one discrete variable $n=t_{1,0}=-t_{2,0}$.
Let us  choose the normalization of the Baker-Akhiezer function with the
help of the condition $\xi_{1,0}= 1$, i.e., in the neighborhood of $P_1$ the \BA function has the form:
\beq
\psi (n,t_{\a,i>0},Q) =
k_1^n\exp \biggl(\sum_{i=1}^{\infty} t_{1,i} k_{1}^{i} \biggr)
\biggl(1+ \sum_{s=1}^{\infty} \xi_{1,s}(n,t) k_1^{-s}\biggr),\label{2.82d1}
\eeq
and in the neighborhood of $P_2$
\beq
\psi (n,t_{\a,i>0},Q) =
k_2^{-n}\exp \biggl(\sum_{i=1}^{\infty} t_{2,i} k_{2}^{i} \biggr)
\biggl(\sum_{s=0}^{\infty} \xi_{2,s}(n,t) k_1^{-s}\biggr),\label{2.82d2}
\eeq
According to Theorem~\ref{th2.1}, the function $\psi$ satisfies two sets of differential equations. The compatibility
conditions (\ref{2.6}) within the each set can be regarded as two copies of the KP hierarchies. In addition
the two-point \BA function satisfies differential difference equation (\ref{2.2d}). The first two of them have the form
\beq\label{nov23}
(\p_{1,1}-T+u)\psi=0,\ \ \ (\p_{2,1}-w T^{-1})\psi=0,
\eeq
where
\beq\label{nov231}
u=(T-1)\xi_{1,1}(n,t) , \ \ \ w=e^{\phi_n-\phi_{n-1}},\ \ e^{\phi_n(t)}=\xi_{2,0}(n,t)
\eeq
The compatibility condition of these equations is equivalent to the $2D$ Toda equation (\ref{2DT}) with $\xi=t_{1,1}$ and $\eta=t_{2,1}$. The explicit formula for $\phi_n$ is a direct corollary of the explicit
formula for the \BA function. The normalization of $\psi$ as in (\ref{2.82d1}) defines the coefficient $c$ in
(\ref{2.101})
\beq\label{BA2} 
\psi=\exp\left(n\Omega_{1,0}+\sum t_{\a,i}\Omega_{\a,i}(P)\right)
{\theta(A(P)+nU+\sum U_{\a,i}t_{\a,i}+Z) \,\theta(Z)
\over \theta(nU+\sum U_{\a,i}t_{\a,i}+Z)\, (\theta(A(P)+Z)} ,
\eeq
If we denote $x=0,\ t=t_{1,1}$ and set $t_{1,i>1}=t_{2,i>0}=0$, then up to a constant in $(x,t)$ factor the formula (\ref{BA2}) coincides with (\ref{p}). Expanding $\psi$ at $P_1$ we get the formula for the coefficient $u$ in
in the first linear equation (\ref{nov231}), which coincides with (\ref{u}). That proved ``the only if'' part
of Theorem~\ref{th1.2}.

\subsubsection*{Example 3. Three-point \BA function}

Starting with three-point case, in which the number of discrete variables is $2$, the \BA function satisfies certain linear difference equations (in addition to the differential and the differential-difference equations (\ref{2.3}), (\ref{2.3d})). The origin of these equations is easy to explain. Indeed, if all the continuous variables vanish, $t_{\a,i>0}=0$, then the Baker-Akhiezer function $\psi_{n,m}(P)$, where $n=-t_{1,0}$, $m=-t_{2,0}$, is a meromorphic function having pole of order $n+m$ at $P_3$ and
zeros of order $n$ and $m$ at $P_1$ and $P_2$ respectively, i.e.,
\beq\label{H0}
\psi_{n,m}\in H^0(D+n(P_3-P_1)+m(P_3-P_2)),\ \   D=\g_1+\cdots+\g_g
\eeq
The functions $\psi_{n+1,m},\psi_{n,m+1},\psi_{n,m}$ are all in the linear space $H^0(D+(n+m+1)P_3-nP_1-mP_2)$.
By Riemann-Roch theorem for a generic $D$ the latter space is $2$-dimensional. Hence, these functions are linear dependent, and they can be normalized such the the linear dependence takes the form (\ref{laxdd}).
The theta-functional formula for the \BA function directly implies formulae (\ref{ud}), (\ref{pd}) and proves
``the only if'' part of Theorem~\ref{th1.3}.

For the first glance it seems that everything here is within the framework of classical algebraic-geometry. What might be new brought to this subject by the soliton theory is understanding
that {\it the discrete variables $t_{\a,0}$ can be replaced by continuous ones}. Of course, if
in the formula (\ref{2.101}) the variable $t_{\a,0}$ is not an integer, then $\psi$ is not a single valued function on $\G$. Nevertheless, because the monodromy properties of $\psi$ do not change if the shift of the argument is integer, it satisfied the same type of linear equations with coefficients given by the same type of formulae. It is necessary to emphasize that in such a form the difference equation becomes functional equation.

\noindent{\bf Remark.}
In the four-point case there is three discrete variables $n$, $m$, $l$. In each two of them the \BA function satisfies
a difference equation. Compatibility of these equations is the BDHE equation (\ref{BDHE}).

\section{Dual Baker-Akhiezer function}

The concept of the dual Baker-Akhiezer function $\psi^+(t,P)$ is universal
and is at the heart of Hirota's bilinear form of soliton equations, and
plays an essential role in our proof of Welters' conjecture.
It is necessary to emphasize that, although the concept is universal,
the definition of the dual \BA function depends on a choice of {\it dual\/}
divisor $D^+=\g_1^++\cdots+\g_g^+$. As it will be shown later
the notion of duality between divisors of $\psi$ and $\psi^+$
reflects a choice of one of the variables $t_{\a,0}$ or $t_{\a,1}$.
In all the cases the pole divisor $D^+$ of the dual \BA function is
defined by the equation
\beq\label{D+D}
D+D^+=K+\kappa\in J(\G)
\eeq
where $K$ is a canonical class and $\kappa$ is a certain degree~$2$
divisor, that encodes the type of duality. Depending on its choice, the
dual \BA function is then defined by the following analytic properties:

i) the function $\psi^+$ (as a function of the variable $P\in \G$) is meromorphic everywhere except for the points $P_{\a}$ and has at most simple poles at the points $\gamma_1^+,\ldots,\gamma_g^+$ (if all
of them are distinct);

(ii) in a neighborhood of the point $P_{\a}$ the function $\psi$ has the
form
\beq
\psi (t,Q) =k^{-t_{\a,0}} \exp \biggl(\sum_{i=1}^{\infty} -t_{\a ,i} k_{\a}^{i}
\biggr) \biggl( \sum_{s=0}^{\infty} \xi_{\a,s}^+(t) k_{\a}^{-s} \biggr),\ \
k_{\a}=k_{\a} (Q) . \label{2.1+}
\eeq
In fact it is the same \BA type function and, therefore, admits the same type of explicit theta-function formula:
\beq
\psi^+(t,P)=c^+(t)\exp\left(-\sum t_{\a,i}\Omega_{\a,i}(P)\right)
{\theta(A(P)-\sum U_{\a,i}t_{\a,i}-Z+\wh \kappa)
\over \theta(A(P)-Z+\wh \kappa)} .\label{2.101+}
\eeq
The basic type of duality and their meaning are explained below in two examples.

\subsubsection*{Example 1. One-point case. Duality for a continuous variable}

The notion of dual \BA function in the one point case was first introduced in \cite{cherednik}.
In this dual divisor is defined by (\ref{D+D}) where $\kappa=2P_1$.
In other words, for a generic effective degree
$g$ divisor $D$ there exists a unique meromorphic differential $d\Omega$ with pole of degree $2$ at $P_1$,
$d\Omega=d(k_1+O(1))$ having zeros at the points $\g_s$; in addition it has $g$ more zeros that are denoted by
$\g_1^+,\ldots,\g_g^+$.

The functions $\psi$ and $\psi^+s$ have essential singularities, their product or products of their derivatives
are meromorphic functions on $\G$. Moreover, from the definition of the duality it follows that after multiplication by corresponding differential $d\Omega$ one gets a meromorphic differential on $\G$ with the only pole at $P_1$. That proves the following statement.

\begin{lem} Let $\psi$ and $\psi^+$ be the Baker-Akhiezer function and its dual. Then the following equations hold:
\beq\label{resba}
\res_{P_1} \left(\psi^+(\p_x^j \psi)\right)d\Omega=0, \ \  j=0,1,\ldots.
\eeq
\end{lem}
Equations (\ref{resba}) allows to express the coefficients $\xi_s^+$ of the expansion of the dual function
$\psi^+$ at $P_1$ as universal differential polynomials in terms of the coefficients $\xi_{s'}$ of the Baker-Akhiezer function. The first such equation is $\xi_1+\xi_1^+=0$.
Another corollary of (\ref{resba}) is infinite number of
bilinear identities for the theta-function, that one obtains after substitution of (\ref{2.101}), (\ref{2.101+})
into (\ref{resba}). These identities are usually called {\it Hirota's bilinear equations}.

\begin{cor} Let $\psi$ be the \BA function and $L_i$ be the linear operator of the form (\ref{2.9}) such that
$(\p_n-L_n)\psi=0$. Then the dual \BA function is a solution of the formal adjoint equation
\beq\label{2.9adj}
\psi^+(\p_n-L_n)=0
\eeq
\end{cor}
Recall that the {\it right action\/} of a differential operator is defined as a formal adjoint action,
i.e., $f^+\p_i=-\p_if^+$ (and the left-hand side of this formula should not
be confused with the more common differentiation-followed-by-multiplication
construction for a differential operator). The proof
of the corollary will be given in the next section.

\subsubsection*{Example 2. Two-point case. Duality for a discrete variable}

In the two-point case, in which there is one discrete variable $n$, the dual divisor $D^+$ is defined by (\ref{D+D}) with $\kappa=P_1+P_2$, i.e., $\g_s$ and $\g_{s'}$ are zeros of a differential $d\Omega$ having
simple poles at the marked points $P_1$ and $P_2$. Without loss of generality we may assume that at these points it has residues $\mp 1$.

\begin{lem}\label{lem3.1} Let $\psi$ and $\psi^+$ be the Baker-Akhiezer function and its dual. Then the following
equations hold:
\beq\label{resbad}
\res_{P_1} \left(\psi^+(T^i \psi)\right)d\Omega=0, \ \  i=1,2\ldots.
\eeq
\end{lem}
By definition of the duality, the differential on the left-hand side of (\ref{resbad}) has pole only at $P_1$.
Hence its residue vanishes. Note also that the differential $\psi^+\psi d\Omega$ has poles at $P_1$ and $P_2$. The constant $c^+$ in the normalization of the dual Baker-Akhiezer function is chosen such that
\beq\label{resba0}
\res_{P_1} \left(\psi^+\psi\right)d\Omega=1.
\eeq
\begin{cor}Let $\psi$ be the \BA function and let $\wh L_n$ be the linear operator of the form
\beq\label{LN}
\wh L_i=T^i+\sum_{j=0}^{n-1}v_j^{(n)}T^j
\eeq
such that   $(\p_{1,i}-\wh L_n)\psi=0$. Then the dual \BA function is a solution of the formal adjoint equation
\beq\label{2.9adj-2Toda}
\psi^+(\p_{1,i}-\wh L_i)=0
\eeq
\end{cor}
As in the case of differential operators,
here and below the right action of a difference operator is defined as formal adjoint action, i.e., $f^+T=T^{-1}f^+$.

\section{Integrable hierarchies}

In its original form equations (\ref{2.6}), (\ref{2.6d}) is just an infinite system of partial differential equation for an infinite number of coefficients of all the operators, depending on infinite number of independent variables called ``times''. Of course, restricting to a finite number of variables one gets an equation or a finite number of equations for a finite number of variables. Some of them are fundamental equations of mathematical physics, and as such deserve special interest. That is true for all three basic equations mentioned above, that is KP, 2D Toda and BDHE. Our next goal is to present the hierarchies of these equations in the form of commuting flows on a certain ``phase spaces'' that are spaces of pseudo\-differential or pseudo\-difference operators.  This form is due to Sato and his coauthors \cite{sato}.

\subsection*{KP hierarchy}

Let $\O$ be a linear space of a formal pseudo\-differential operators in the variable $x$, i.e., formal series
\beq\label{pseudo}
\D=\sum_{s=-N}^{\infty} v_s(x)\p_x^{-s}
\eeq
By definition the coefficient $v_1$ at $\p_x^{-1}$ in (\ref{pseudo}) is called the residue of $\D$
\beq\label{Dres}
v_1:=\res_\p \D.
\eeq
The commutator relations $\p_x\cdot v(x)=v_x(x)+v(x)\p_x$ and $\p_x^{-1}\cdot v(x)=v(x)\p_x^{-1}-v_x(x)\p_x^{-2}+v_xx(x)\p_x^{-2}$ define on $\O$ a structure of associative ring. For any
pseudo\-differential operator $\D$ its differential part is defined as the unique differential operator such that
$\D-\D_+=\D_-=O(\p_x^{-1})$, i.e., for $\D$ as in (\ref{pseudo}) its differential part is equal to
\beq\label{pseudo+}
\D_+=\sum_{s=-N}^{0} v_s(x)\p_x^{-s}
\eeq
The KP hierarchy is defined on the space $\P$ of monic pseudo\-differential operators of order $1$, i.e., of the operators of the form
\beq\label{Lkp}
\L=\p_x+\sum_{s=1}^\infty v_s(x)\p_x^{-s}
\eeq
\begin{prop} The equations
\beq\label{kphier}
\p_{i} \L=[\L^i_+,\L]
\eeq
define commuting flows on the space $\P$.
\end{prop}
{\it Proof.}  The left-hand side of equation (\ref{kphier}) is a pseudo\-differential operator $\p_i \L=\sum_{s\ge1} (\p_i v_s)\p^{-s}$ of order at most $-1$. Therefore, (\ref{kphier}) is well-defined if and only if the right-hand side is a pseudo\-differential operator of order at most $-1$. To show this, notice, that the identity $[\L^i,\L]=0$ implies $[\L^i_+,\L]=-[\L^i_-,\L]$. Be definition $\L^i_-$ is an operator of order at most $-1$. Hence, $[\L^i_-,\L]$ is also of order at most $-1$.

For the proof of the second statement of the proposition it is necessary to show that equations (\ref{kphier})
imply the equation
\beq\label{kphier1}
[\p_i-\L^i_+,\p_j-\L^j_+]=\p_i \L^j_+-\p_j \L^i_++[\L^j_+,\L^i_+]=0
\eeq
The left-hand side of (\ref{kphier1}) is a differential operator. Therefore, in order to show that it vanish, it is enough to show that it is a pseudo\-differential operator of order at most $-1$. From (\ref{kphier})
it follows that $\p_i \L^j=[\L_+^i,\L^j]$ Then using the the identity $[\L^i,\L^j]=0$ we have
\beq\label{kphier2}
\p_i \L^j_+=[\L_+^i,\L^j]-\p_i \L^j_-=[\L^j,\L_-^i]+O(\p_x^{-1})=
[\L^j_+,\L_-^i]+O(\p_x^{-1})
\eeq
Similarly,
\beq\label{kphier3}
[\L^i_+,\L^j_+]=[\L^j_+,\L^i_-]-[\L^j_+,\L^i_-]+O(\p_x^{-1})
\eeq
Substituting (\ref{kphier2}), (\ref{kphier3}) into (\ref{kphier1}) completes the proof of the proposition.

The operator $\L^2_+$ has the form $\p_x^2-u(x,y)$, with $u=-2v_1$ where $v_1$ is the coefficient at $\p_x^{-1}$ of $\L$, i.e., $v_1=\res_\p\L$. Equations (\ref{kphier1}) with $j=2$ have the form
\beq\label{kphier4}
\p_{t_m} u=[\p_y-\p_x^2+u,\L^m_+]=-[\p_y-\p_x^2+u, \L^m_-]=2\p_x F_m\,,
\eeq
where
$$	F_m:=\res_\p \L^m.
$$

\paragraph*{Important remark} At first glance the system (\ref{kphier4}) looks like a system of commuting
evolution equations, but it is not. The right-hand side of (\ref{kphier4}) are universal
differential polynomials in $v_i$. In general there is no way to reconstruct from one function $u(x,y)$
an infinite set of functions $v_i(x)$ of one variable. It can be done only under ceratin assumptions.
In \cite{kp} that was done in the case when $u(x,y)$ is a periodic function of the variables $x$ and $y$.
To some extend the main part in the proof of the first case of Welter's conjecture can be seen as
the proof of the equivalence of (\ref{kphier}) and (\ref{kphier4}) in the case when $u$ is as in
the statement of Theorem~\ref{th1.1}.

For further use let us present some other basic notations and construction. The first one is the notion of
wave function.
\begin{lem} Let $\L$ be a monic pseudo\-differential operator of the form (\ref{Lkp}). Then the equation $\L\psi=k\psi$ has a unique solution of the form
\beq\label{kphier6}
\psi=e^{kx}\biggl(1+\sum_{s=1}^\infty \xi_s(x)k^{-s}\biggr)
\eeq normalized by the condition $\xi_s(0)=0$.
\end{lem}
The proof is elementary. Substituting (\ref{kphier6}) into the equation gives a system of equations having the form $p_x\xi_s=R_s(v_k,\xi_s')$ with $k,s'<s$. Therefore, they  uniquely define $\xi_s$,if the initial conditions are fixed.

The wave function is then define the wave operator
\beq\label{waveoper}
\Phi=1+\sum_{s=1}^\infty\varphi_s(x)\p_x^{-s}
\eeq
by the equation $\psi=\Phi e^{kx}$.
Notice, that the last equation implies
\beq\label{Lphi}
\L=\Phi\cdot \p_x \cdot \Phi^{-1}
\eeq
The formal dual wave function is given by the formula
\beq\label{psidual}
\psi^+=e^{-kx}\biggl(1+\sum_{s=1}^\infty \xi_s^+(x) k^{-s}\biggr):=e^{-kx}\Phi^{-1}
\eeq
is a solution of the formal adjoint equation $\psi^+\L=k\psi^+$

The defining property of the dual wave function are equations that we proved for the dual \BA function
in the previous section. Namely,
\begin{lem} Let $\psi$ be a wave function and $\psi^+$ its dual. Then the equations
\beq\label{dualnov}
\res_k (\psi^+ (\p_x^n\psi))\,dk =0, \ \ \ n=0,1,\ldots
\eeq
hold.
\end{lem}
The proof is a direct corollary of the identity
\beq\label{dic}
\res_k \left(e^{-kx}\D_1\right)\left(\D_2e^{kx}\right)\,dk=\res_{\p}\left(\D_2\D_1\right),
\eeq
which holds for any pair of pseudo\-differential operators (for details see \cite{sato,dikii}).

In the same way one can show that the product of the wave function and its dual is a generating series for
the right-hand sides of the hierarchy (\ref{kphier4}).
\begin{lem}\label{lem4.3} The coefficients of the expansion
\beq\label{prodnov}
\psi^+\psi=1+\sum_{s=2}^\infty J_sk^{-s}
\eeq
are given by $J_{n+1}=F_n=\res_\p \L^{n}$.
\end{lem}
{\it Proof.} From the definition of $\L$ it follows that
\beq\label{20}
\res_k\left(\psi^+(\L^n\psi)\right)\, dk=\res_k\left(\psi^+k^n\psi\right)\,dk=J_{n+1}.
\eeq
On the other hand, using the identity (\ref{dic}) we get
\beq\label{201}
\res_k(\psi^+\L^n\psi)\,dk=\res_k\left(e^{-kx}\Phi^{-1}\right)\left(\L^n\Phi e^{kx}\right)\,dk=
\res_{\p}\L^n=F_n.
\eeq
The lemma is proved.

\subsection*{2D Toda hierarchy}
In the two-point case there are two sets of continuous variables and one discrete variable which we
denote by $x$. It is instructive enough to consider the hierarchy of equations corresponding to one set of continuous times associated with one marked point. In this subsection we present the definition of the hierarchy of the differential-difference equations (\ref{2.6d}) in the form of the commuting flows on the space $\P$
of the pseudo\-difference operators of the form
\beq\label{LLnov}
\L=T+\sum_{s=0}^{\infty} w_s(x)T^{-s}, \ \ T=e^{\p_x}
\eeq
In the ring of the pseudo\-difference operators
\beq\label{psdif}
\D=\sum_{s=-N}^{\infty} v_s(x) T^{-s}
\eeq
the notion of the residue as follows:
\beq\label{resdif}
\res_T \D:=v_0
\eeq
For any pseudo\-differential operator $\D$ its positive part is defined as
the difference  operator such that $\D_-:=\D-\D_+=O(T^{-1}$, i.e., if $\D$ is as in (\ref{psdif}), then
\beq\label{+def}
\D_+:=\sum_{s=-N}^{-1} w_s(x)T^{-s}
\eeq
\begin{prop} The equations
\beq\label{2dh}
\p_{i} \L=[\L^i_+,\L]
\eeq
define commuting flows on the space $\P$.
\end{prop}
The proof of the first statement goes along the same lines as in the case of KP hierarchy. The proof of the second statement that (\ref{2dh}) implies
\beq\label{2dh1}
[\p_i-\L^i_+,\p_j-\L^j_+]=0
\eeq
is also identical. The first operator $\L^+$ is of the form $\L_+=T-u$ with $u=w_0$. The equation (\ref{2dh}) for $i=1$ gives $\p_t u=-w_1$, where $w_1=\res_T\,\L\,T$. Here and below $t=t_1$.
For further use, let us present the equation
\beq\label{FF1}
\p_t F_m=(1-T)F_m^1,
\eeq
where
$$F_m=\res_T\L^m,\quad F_m^1=\res_T \L^mT,$$
which directly follows from the comparison of residues of two side of the equality $\p_t \L^m=[\L+,\L_m]$.
The commutativity equations (\ref{2dh1}) imply that the evolution of $u$ with respect to all the other times
\beq\label{2dev}
\p_{t_m} u= -(T-1) F_m^1=-\p_t F_m
\eeq
As in the KP case, in general the last equations can not be regarded as well-defined hierarchy on the
space of one function $u(x,t)$ because the definition of $F_m$ involves other coefficients of $\L$. The main part of the proof of the second case of Welters' conjecture can be seen as a reconstruction of $\L$ in terms of $u$ under the assumption of Theorem~\ref{th1.2}.

We conclude this section by providing a necessary definitions and identities, which are just discrete analog of that above. Namely, the wave function is a solution of the equation $\L\psi=k\psi$ of the
form
\beq\label{psidis}
\psi=k^x\biggl(1+\sum_s \xi_s(x) k^{-s}\biggr)
\eeq
It defines a unique wave operator by the equation
\beq\label{Snov}
\psi=\Phi k^x,\ \ \Phi=1+\sum_{s=1}^{\infty}\f_s(x)T^{-s}.
\eeq
Then, the dual wave function is defined by the left action of the operator $\Phi^{-1}$:
$\psi^+=k^{-x}\Phi^{-1}$.
Recall that the left action of a pseudo\-difference operator is the formal adjoint action
under which the left action of $T$ on a function $f$ is $(fT)=T^{-1}f$.
\begin{lem} The coefficient of the product
\beq\label{J-2Toda}
\psi^+\psi=1+\sum_{s=1}^{\infty}J_s(Z,t)\,k^{-s}
\eeq
are equal to $J_n=F_n=\res_T \L^n$.
\end{lem}
{\it Proof.}
 From the definition of $\L$ it follows that
\beq\label{20-2Toda}
\res_k\left(\psi^+(\L^n\psi)\right)k^{-1}dk=\res_k\left(\psi^+k^n\psi\right)k^{-1}dk
=J_{n}.
\eeq
On the other hand, using the identity
\beq\label{dic-2Toda}
\res_k \left(k^{-x}\D_1\right)\left(\D_2k^x\right)k^{-1}dk=
\res_T\left(\D_2\D_1\right)
\eeq
which is the 2D Toda analogue of (\ref{dic}),
we get
\beq\label{201-2Toda}
\res_k(\psi^+\L^n\psi)k^{-1}dk=
\res_k\left(k^{-x}\Phi^{-1}\right)\left(\L^n\Phi k^x\right)k^{-1}dk=
\res_T\L^n=F_n.
\eeq
Therefore, $F_n=J_{n}$ and the lemma is proved.

\section{Commuting differential and difference operators. }

In the previous section hierarchies of the KP and 2D Toda equations were defined as systems of commuting flows
on the spaces of pseudo\-differential or pseudo\-difference operators, respectively. Consider now the subspace $\O_n\subset \O$ of operators whose $n$-th power is a differential (difference) operator $L_n$, i.e., $\L^n=L_n$
or equivalently $\L^n_-=0$. The latter directly implies that $\p_{t_n} \L=0$. In other words the subspace $\O_n$ is the subspace of stationary points of the $n$-th flow of the hierarchy. It has finite functional dimension and can be simply identified with the space of all monic differential (difference) operators because any such operator
$L_n$ uniquely defines the corresponding pseudo\-differential $\L=L_n^{1/n}$.
The subspace $\O_n$ is invariant with respect to all the other flows. Their restriction on $\O_n$ is a closed
system of evolution equations on a space of finite-number of unknown functions and can be represented in the form
$\p_i L_n=[L^{i/n}_{n,+}, L_n]$. For $n=2$ the corresponding reduction of the KP hierarchy is equivalent to the
hierarchy of the KdV equation $4u_t=6uu_x+u_{xxx}$. An attempt to find explicit periodic solutions of the KdV equation had led Novikov in to the idea to consider further reduction to stationary points of one of the ``higher'' KdV flows. In terms of the original KP hierarchy that is a subspace stationary for two flows
of the hierarchy (or two linear combinations of basic flows). The corresponding subspace is
the space of differential order $n$ monic ordinary differential operator $L_n$ such that there exists operator $L_m$ commuting with $L_n$ of order $m$ (not multiple of $n$), i.e., the space of solutions of a system (\ref{comoper}). As it was mentioned in the introduction, the problem of classification of commuting ordinary differential operators as pure algebraic problem was consider
in remarkable works by Burchnall and Chaundy \cite{ch}.

Briefly the key points of their proof of the statement that a pair of such operators is always satisfy algebraic relation
\beq
R(L_n,L_m)=0. \label{12}
\eeq
are the following. The commutativity of $L_n$ and $L_m$ implies that the space
${V}(\lambda)$ of solutions of the ordinary linear equation $L_n y(x)=\lambda y(x)$
is invariant with respect to the operator $L_m$.  The matrix elements
$L_m^{ij}$ of the corresponding finite dimensional linear operator $L_m(\lambda)$
\beq
L_m|_{V(\lambda)} =L_m(\lambda)\colon V(\lambda)\longmapsto V(\lambda) \label{15}
\eeq
in the canonical basis
$
c_i(x,\lambda,x_0)\in \L(\lambda), \ \ \
c_i(x,\lambda,x_0)|_{x=x_0}=\delta_{ij}, $
are polynomial functions in the variable $\lambda$. They depend on the choice
of the normalization point $x=x_0$, i.e., $L_m^{ij}=L_m^{ij}(\lambda, x_0)$.
The characteristic polynomial
\beq
R(\lambda, \mu)=\det(\mu-L_m^{ij}(\lambda, x_0))\label{17}
\eeq
is a polynomial in both variables $\lambda$ and $\mu$ and does not depend on
$x_0$.

According to the property of characteristic polynomials we have
$$
R(L_n,L_m)y(x,\lambda)=0.
$$
Notice, that
$R(L_n,L_m)$ is an ordinary differential operator. Therefore, if it is
not equal to zero then its kernel is finite dimensional. Hence,
the last equation valid for all $\lambda$ implies (\ref{12}), and
the first statement of \cite{ch} is proved.

The equation $R(\lambda,\mu)=0$ defines affine part of an algebraic curve. Let us show that it is always compactified by one {\it smooth} point $P_0$. Indeed the equation $L_n\psi=k^n \psi$ has always a unique
formal wave solution, i.e., a solution of the form (\ref{kphier6}) normalized by the conditions $\xi_s(0)$=0.
Moreover, any solution of the latter equation of the form $e^{kx}\cdot(\hbox{Laurent series in }k^{-1})$ is equal to $\psi(x,k)c(k)$, where $c(k)$ is a constant Laurent series. The operator $L_m$ commutes with $L_n$, therefore $L_m\psi$ is also a solution to the same equation. Hence, there exists a Laurent series
\beq\label{anov}
a_m(k)=k^m+\sum_{s=-m+1}^\infty a_{m,s} k^{-s}
\eeq
such that $\L_m\psi=a(k)\psi(x,k)$, i.e., $\psi$ is a formal common eigenfunction of the operators $L_n,L_m$.
That implies the following expansion of the characteristic equation at infinity $\lambda\to infty$:
\beq\label{rexpans}
R(\lambda,\mu)=\prod_{i=0}^{n-1}(\mu-a(k_i)),\ \  k_i^n=\lambda.
\eeq
Now we are ready to explain a role of the condition under which Burchnall and Chaundy where able to make the next step. Namely, the condition that orders of operators are co-prime. The leading coefficient of $a(k)$
is $k^m$. Hence, if $(n,m)=1$ then in the neighborhood of the infinite (and, therefore, almost
everywhere else) the operator $L_n(\lambda)$ has $n$-distinct eigenvalues, and is diagonalizable, i.e.,
for each generic point $P=(\lambda,\mu)\in \G$ there is a unique eigenfunction $\psi(x,P;x_0)$ of the
operators $L_n,L_m$ normalized by the condition $\psi(x_0,P;x_0)=1$. It can be written as
\beq\label{psic}
\psi(x,P;x_0)=\sum_{i=0}^{n-1} h_i(P,x_0)c_i(x,\lambda;x_0), \ \ h_0(P,x_0)=1,
\eeq
where $c_i$ are canonical basis of solution to the equation $L_ny=\lambda y$ defined above and $h_i$ are
coordinates of the eigenvector of the matrix $\L_m(\lambda)$. They are rational expressions in $\lambda$ and $\mu$, and, therefore are meromorphic functions of $P\in \G$ (if $\G$ is smooth, otherwise they become meromorphic on an normalization of $\G$). The functions $c_i$, as solutions of the initial value problem, are entire function of the variable $\lambda$. Hence, $\psi$ in an affine part of $\G$ is a meromorphic function
with poles that are independent of $x$ (but depend on the normalization point $x=x_0$). If $\G$ is smooth than
their number is equal to the genus $g$ of $\G$. By definition of the canonical basis we have that
$\psi_x (x,P)\psi^{-1}(x,P)|_{x=x_0}=h_1(P,x_0)$. The asymptotic of $h_1$ can be easy found using the formal
wave solution. It equals $h_1=k+(O(k^{-1}))$. Therefore  $\psi=\exp \left(\int_{x_0} h_1(x,P)dx\right)$
has at $P_0$ exponential singularity and is a Baker-Akhiezer function (with the shift of $x$ by $x_0$).

\begin{theo} \cite{ch,kr1,kr2,mum} There is a natural correspondence
\beq\label{corr}
\A\longleftrightarrow \{\G,P_0, [k^{-1}]_1, \F\}
\eeq
between {\it regular} at $x=0$ commutative rings $\A$ of ordinary linear
differential operators containing a pair of monic operators of co-prime orders, and
sets of algebraic-geometrical data $\{\G,P_0, [k^{-1}]_1, \F\}$, where $\G$ is an
algebraic curve with a fixed
first jet $[k^{-1}]_1$ of a local coordinate $k^{-1}$ in the neighborhood of a smooth
point $P_0\in\G$ and $\F$ is a torsion-free rank 1 sheaf on $\G$ such that
\beq\label{sheaf}
H^0(\G,\F)=H^1(\G,\F)=0.
\eeq
The correspondence becomes one-to-one if the rings $\A$ are considered modulo conjugation
$\A'=g(x)\A g^{-1}(x)$.
\end{theo}
Note that in \cite{kr1,kr2,ch} the main attention was paid to the generic case of
the commutative rings corresponding to smooth algebraic curves.
The invariant formulation of the correspondence given above is due to Mumford \cite{mum}.

The algebraic curve $\G$ is called the spectral curve of $\A$.
The ring $\A$ is isomorphic to the ring $A(\G,P_0)$ of meromorphic functions
on $\G$ with the only pole at the point $P_0$. The isomorphism is defined by
the equation
\beq\label{z2}
L_a\psi_0=a\psi_0, \ \ L_a\in \A, \ a\in A(\G,P_0).
\eeq
Here $\psi_0$ is a common eigenfunction of the commuting operators. At $x=0$ it is
a section of the sheaf $\F\otimes\O(-P_0)$.

\noindent{\bf Remark.} As we have seen above, the construction of the correspondence (\ref{corr})
depends on a choice of initial point $x_0=0$. The spectral curve and the sheaf $\F$
are defined by the evaluations of the coefficients of generators of $\A$ and a finite
number of their derivatives at the initial point. In fact, the spectral curve
is independent on the choice of $x_0$, but the sheaf does depend on it, i.e., $\F=\F_{x_0}$.

Using the shift of the initial point it is easy to show that the correspondence
(\ref{corr}) extends to the commutative rings of operators whose coefficients are
{\it meromorphic} functions of $x$ at $x=0$. The rings of operators having poles at $x=0$
correspond to sheaves for which the condition (\ref{sheaf}) is violated.

\noindent{\bf Remark.} In their original paper Burchnall and Chaundy stressed that there is no approach to a classification of commutative differential operators whose ordered are not co-prime. The classification of
commutative rings of ordinary differential operators was completed in \cite{kr-com}, where it was shown that
a maximal ring $\A$ of commuting differential operators is uniquely defined by an algebraic curve with marked point,
the first jet of local coordinate at the marked point, and if the curve is smooth by the rank $k$ and degree $rg$
vector bundle. In addition it depends on $r-1$ arbitrary functions of one variable. Here $k$ is the rank of $\A$
defined as the greatest common divisor of the orders of commuting operators.

\subsection*{Commuting difference operators}

A theory of commuting difference operators containing a pair of operators of co-prime orders
was developed in  \cite{mum,kr-dif}. It is analogous to the theory of rank 1 commuting
(Relatively recently this theory was generalized to the case of commuting difference operators of arbitrary rank
in \cite{n-kr}.) For further use we present here the classification of commutative differential operators
of the form
\beq\label{lform}
L_n=T^n+\sum_{s=1}^{n-1} u_i(x)T^i
\eeq

\begin{theo} (\cite{mum,kr-dif}) Let $\A$ be a maximum commutative ring of ordinary difference
operators of the form (\ref{lform}) containing a pair of operators of co-prime orderes. Then there is an irreducible algebraic curve $\G$, such that  the ring $\A^Z$ is isomorphic to the ring $A(\G,P_+,P_-)$ of the meromorphic functions on $\G$ with the only pole at a smooth point $P_+$, vanishing at another smooth point
$P_-$. The ring is uniquely defined by a torsion-free rank 1 sheaves $\F$ on $\G$
such that
\beq\label{sheaf-2Toda}
h^0(\G,\F(nP_+-nP_-))=h^1(\G,\F(nP_+-nP_-))=0.
\eeq
The correspondence becomes one-to-one if the rings $\A$ are considered modulo conjugation
$\A'=g(x)\A g^{-1}(x)$.
\end{theo}
{\bf Remark.}
As in the continuous case the construction of the correspondence
depends on a choice of initial point $x_0=0$. The spectral curve and the sheaf $\F$
are defined by the evaluations of the coefficients of generators of $\A$
at a finite number of points of the form $x_0+n$.
In fact, the spectral curve is independent on the choice of $x_0$, but the sheaf
does depend on it, i.e., $\F=\F_{x_0}$.

Using the shift of the initial point it is easy to show that the correspondence
(\ref{corr}) extends to the commutative rings of operators whose coefficients are
{\it meromorphic} functions of $x$. The rings of operators having poles at $x=0$
correspond to sheaves for which the condition (\ref{sheaf-2Toda}) for $n=0$
is violated.
\section{Proof of Welters' conjecture}

As it was mentioned in the introduction the proof of all the particular cases of Welters' trisecant conjecture
uses different hierarchies: the KP, the 2D Toda, and BDHE. In each case there are some specific difficulties but
the main ideas and structures of the proof are the same. In all the cases the first step is to construct
the wave solution. It is necessary to emphasize that it is not a wave solution to the ordinary
pseudo\-differential or pseudo\-difference operators discussed in Section 4. The corresponding wave solutions
are defined as formal solutions to a {\it partial differential equation}. In this case there is no way
to define such a solution in a unique way without additional assumption on a global structure of the
coefficients of the equation. As an instructive example we present in this section the proof of the first
particular case of Welters' conjecture, namely, the proof of Theorem~\ref{th1.1}.

First, we prove the implication $(A)\to (C)$.
Let $\tau(x,y)$ be a holomorphic function of the variable $x$ in some open domain
$D\in \mathbb C$ smoothly depending on a parameter $y$. Suppose that for each $y$
the zeros of $\tau$ are simple,
\beq\label{xi}
\tau(x_i(y),y)=0,\quad\tau_x(x_i(y),y)\neq 0.
\eeq
\begin{lem} (\cite{flex})
If equation (\ref{lax0}) with the potential
$u=-2\p_x^2\ln \tau(x,y)$ has a meromorphic in $D$ solution $\psi_0(x,y)$, then
equations (\ref{cm50}) hold.
\end{lem}
{\it Proof.}
Consider the Laurent expansions of $\psi_0$ and $u$ in the neighborhood of one of the zeros
$x_i$ of $\tau$:
\beq\label{ue}
\begin{split}
u&={2\over (x-x_i)^2}+v_i+w_i(x-x_i)+\ldots;\\
\psi_0&={\a_i\over x-x_i}+\b_i+\g_i(x-x_i)+\delta_i(x-x_i)^2+\ldots.
\end{split}
\eeq
(All coefficients in these expansions are smooth functions of the variable $y$).
Substitution of (\ref{ue}) in (\ref{lax0}) gives a system of equations. The first three of them
are
\beq\label{eq1}
\a_i \dot x_i+2\b_i=0;\
\dot\a_i+\a_i v_i+2\g_i=0;\
\dot\b_i+v_i\b_i-\g_i\dot x_i+\a_i w_i=0.
\eeq
Taking the $y$-derivative of the first equation and using two others we get
(\ref{cm50}).

Let us show that equations (\ref{cm50}) are sufficient for the existence of
meromorphic wave solutions, i.e., solutions of the form (\ref{ps}).
\begin{lem}
Suppose that equations (\ref{cm50}) for the zeros of $\tau(x,y)$ hold. Then there exist
meromorphic wave solutions of equation (\ref{lax0}) that have simple poles at $x_i$ and
are holomorphic everywhere else.
\end{lem}
{\it Proof.} Substitution of (\ref{ps}) into (\ref{lax0}) gives a recurrent system of
equations
\beq\label{xis}
2\xi_{s+1}'=\p_y\xi_s+u\xi_s-\xi_s''
\eeq
We are going to prove by induction that this system has meromorphic solutions with
simple poles at all the zeros $x_i$ of $\tau$.

Let us expand $\xi_s$ at $x_i$:
\beq\label{5}
\xi_s={r_s\over x-x_i}+r_{s0}+r_{s1}(x-x_i)\,,
\eeq
where for brevity we omit the index $i$ in the notations for the coefficients of this
expansion.
Suppose that $\xi_s$ are defined and equation (\ref{xis}) has a meromorphic solution.
Then the right-hand side of (\ref{xis}) has the zero residue at $x=x_i$, i.e.,
\beq\label{res}
{\rm res}_{x_i}\left(\p_y\xi_s+u\xi_s-\xi_s''\right)=\dot r_s+v_ir_s+2r_{s1}=0
\eeq
We need to show that the residue of the next equation vanishes also.
 From (\ref{xis}) it follows that the coefficients of the Laurent expansion for $\xi_{s+1}$
are equal to
\beq\label{6}
r_{s+1}=-\dot x_ir_s-2r_{s0},
\eeq
\beq\label{7}
 2r_{s+1,1}=\dot r_{s0}-r_{s1}+w_ir_s+v_ir_{s0}\,.
\eeq
These equations imply
\beq
\dot r_{s+1}+v_ir_{s+1}+2r_{s+1,1}=-r_s(\ddot x_i-2w_i)-\dot x_i(\dot r_s-v_ir_ss+2r_{s1})=0,
\eeq
and the lemma is proved.

\subsection*{$\lambda$-periodic wave solutions}

Our next goal is to fix a {\it translation-invariant} normalization of $\xi_s$
which defines wave functions uniquely up to a $x$-independent factor.
It is instructive to consider first the case of the periodic potentials $u(x+1,y)=u(x,y)$
(see details in \cite{kp}).

Equations (\ref{xis}) are solved recursively by the formulae
\beq
\xi_{s+1}(x,y)=c_{s+1}(y)+\xi_{s+1}^0(x,y)\,,\label{kp1}
\eeq
\beq\label{kp2}
\xi_{s+1}^0(x,y)={1\over 2}\int_{x_0}^x (\p_y \xi_s-\xi_s''+u\xi_s)\,dx\, ,
\eeq
where $c_s(y)$ are {\it arbitrary} functions of the variable $y$.
Let us show that the periodicity condition $\xi_s(x+1,y)=\xi_s(x,y)$
defines the functions $c_s(y)$ uniquely up to an additive constant.
Assume that $\xi_{s-1}$ is known  and satisfies the condition that the corresponding
function $\xi_s^0$ is periodic.
The choice of the function $c_s(y)$ does not affect the periodicity property of
$\xi_s$, but it does affect the periodicity in $x$ of the function
$\xi_{s+1}^0(x,y)$. In order to make  $\xi_{s+1}^0(x,y)$ periodic,
the function $c_s(y)$ should satisfy the linear differential equation
\beq\label{kp4}
\p_y c_s(y)+B(y)\,c_s(y)+\int_{x_0}^{x_0+1} \left(\p_y \xi_s^0(x,y)+
u(x,y)\,\xi_s^0(x,y)\right)\,dx\ ,
\eeq
where $B(y)=\int_{x_0}^{x_0+1} u\, dx$.
This defines $c_s$ uniquely up to a constant.

In the general case, when $u$ is quasi-periodic, the normalization of the wave functions
is defined along the same lines.

Let $Y_U=\langle \bC U\rangle$ be the Zariski closure of the group
$\bC U=\{Ux\mid x\in\bC\}$ in $X$. Shifting $Y_U$ if needed, we may assume, without loss of generality, that
$Y_U$ is not in the singular locus, $Y_U\not\subset\Sigma$. Then, for a sufficiently small
$y$, we have $Y_U+Vy\notin\Sigma$ as well.
Consider the restriction of the theta-function onto the affine subspace
$\mathbb C^d+Vy$, where
$\mathbb C^d:=($the identity component of $\pi^{-1}(Y_U))$, and
$\pi\colon \mathbb C^g\to X=\mathbb C^g/\Lambda$ is the universal covering map of $X$:
\beq\label{ttt1}
\tau (z,y)=\theta(z+Vy), \ \ z\in \mathbb C^d.
\eeq
The function $u(z,y)=-2\partial_1^2\ln \tau$ is periodic with respect to the lattice
$\Lambda_U=\Lambda\cap \mathbb C^d$ and, for fixed $y$, has a double pole along the divisor
$\Theta^{\,U}(y)=\left(\Theta-Vy\right)\cap \mathbb C^d$.

\begin{lem}\label{lem6.3} Let equations (\ref{cm50}) for zeros of $\tau(Ux+z,y)$ hold and let
$\l$ be a vector of the sublattice $\Lambda_U=\Lambda\cap \mathbb C^d\subset \mathbb C^g$. Then:

(i) equation (\ref{lax0}) with the potential $u(Ux+z,y)$
has a wave solution of  the form $\psi=e^{kx+k^2y}\phi(Ux+z,y,k)$
such that the coefficients $\xi_s(z,y)$ of the formal series
\beq\label{psi2}
\phi(z,y,k)=e^{by}\biggl(1+\sum_{s=1}^{\infty}\xi_s(z,y)\, k^{-s}\biggr)
\eeq
are $\l$-periodic meromorphic functions of the variable $z\in \mathbb C^d$
with a simple pole at
the divisor $\Theta^U(y)$,
\beq\label{v1}
\xi_s(z+\l,y)=\xi_s(z,y)={\tau_s(z,y)\over \tau(z,y)}\, ;
\eeq

(ii) $\phi(z,y,k)$ is unique up to a factor $\rho(z,k)$ that is $\partial_U$-invariant
and holomorphic in $z$,
\beq\label{v2}
\phi_1(z,y,k)=\phi(z,y,k)\rho(z,k), \ \partial_U\rho=0.
\eeq
\end{lem}
{\it Proof.} The functions $\xi_s(z)$ are defined recursively by the equations
\beq\label{xis1}
2\partial_U\xi_{s+1}=\p_y\xi_s+(u+b)\xi_s-\partial_U^2\xi_s.
\eeq
A particular solution of the first equation $2\partial_U\xi_1=u+b$ is given by the formula
\beq\label{v5}
2\xi_1^0=-2\partial_U\ln \tau +(l,z)\ b,
\eeq
where $(l,z)$ is a linear form on $\mathbb C^d$ given by the scalar product of $z$ with
a vector $l\in \mathbb C^d$ such that $(l,U)=1$. By definition, the vector $\l$
is in $Y_U$. Therefore, $(l,\l)\neq 0$. The periodicity
condition for $\xi_1^0$ defines the constant $b$
\beq\label{v6}
(l,\lambda)b=(2\partial_U\ln \tau(z+\l,y)-2\partial_U\ln \tau(z,y))\,,
\eeq
which depends only on a choice of the lattice vector $\l$.
A change of the potential by an additive constant does not affect the results of the
previous lemma. Therefore, equations (\ref{cm50}) are sufficient for the local solvability
of (\ref{xis1}) in any domain, where $\tau(z+Ux,y)$ has simple zeros, i.e., outside
of the set $\Theta_1^{\,U}(y)=\left(\Theta_1-Vy\right)\cap \mathbb C^d$, where
$\Theta_1=\Theta\cap \partial_U\Theta$.
This set  does not contain a $\partial_U$-invariant line because any such line is dense in
$Y_U$. Therefore, the sheaf $\V_0$ of $\partial_U$-invariant meromorphic functions on
$\mathbb{C}^d\setminus \Theta^{\,U}_1(y)$ with poles along the divisor $\Theta^{\,U}(y)$
coincides with the sheaf of holomorphic $\partial_U$-invariant functions.
That implies the vanishing of $H^1({C}^d\setminus \Theta^{\,U}_1(y), \V_0)$ and the
existence of global meromorphic solutions $\xi_s^0$ of (\ref{xis1}) which have a simple
pole at the divisor $\Theta^{\,U}(y)$ (see details in \cite{arbarello,shiota}).
If $\xi_s^0$ are fixed, then the general global meromorphic solutions are given by the
formula $\xi_s=\xi_s^0+c_s$, where the constant of integration $c_s(z,y)$
is a holomorphic $\p_U$-invariant function of the variable $z$.

Let us assume, as in the example above, that a $\l$-periodic solution
$\xi_{s-1}$ is known  and that it satisfies the condition that there exists a periodic
solution $\xi_s^0$ of the next equation. Let $\xi_{s+1}^*$ be a solution of
(\ref{xis1}) for fixed $\xi_s^0$. Then it is easy to see that the function
\beq\label{v7}
\xi_{s+1}^0(z,y)=\xi_{s+1}^*(z,y)+c_s(z,y)\,\xi_1^0(z,y)+{(l,z)\over 2}\p_y c_s (z,y),
\eeq
is a solution of (\ref{xis1}) for $\xi_s=\xi_s^0+c_s$.
A choice of a $\l$-periodic $\p_U$-invariant function $c_s(z,y)$ does not affect the
periodicity property of $\xi_s$, but it does affect the periodicity of the function
$\xi_{s+1}^0$. In order to make  $\xi_{s+1}^0$ periodic,
the function $c_s(z,y)$ should satisfy the linear differential equation
\beq\label{kp41}
(l,\lambda)\p_y c_s(z,y)=2\xi_{s+1}^*(z+\l,y)-2\xi_{s+1}^*(z,y)\,.
\eeq
This equation, together with an initial condition $c_s(z)=c_s(z,0)$ uniquely defines $c_s(x,y)$.
The induction step is then completed. We have shown that the ratio of two
periodic formal series $\phi_1$ and $\phi$ is $y$-independent. Therefore,
equation (\ref{v2}), where $\rho(z,k)$ is defined by the evaluation of
the both sides at $y=0$, holds. The lemma is thus proven.

\begin{cor} Let $\l_1,\ldots,\l_d$ be a set of linear independent vectors of
the lattice $\Lambda_U$ and let $z_0$ be a point of $\mathbb C^d$. Then, under
the assumptions of the previous lemma,
 there is a unique wave solution of equation (\ref{lax0}) such that
the corresponding formal series $\phi(z,y,k;z_0)$ is quasi-periodic with respect
to $\Lambda_U$, i.e., for
$\l\in \Lambda_U$
\beq\label{v10}
\phi(z+\l,y,k;z_0)=\phi(z,y,k;z_0)\,\mu_{\l}(k)
\eeq
and satisfies the normalization conditions
\beq\label{v11}
\mu_{\l_i}(k)=1,\ \ \ \phi(z_0,0,k;z_0)=1.
\eeq
\end{cor}
The proof is identical to that of the part (b) of the Lemma 12 in \cite{shiota}.
Let us briefly present its main steps. As shown above, there exist wave solutions
corresponding to $\phi$ which are $\l_1$-periodic. Moreover, from the statement (ii) above
it follows that for any $\l'\in \Lambda_U$
\beq\label{v12}
\phi(z+\l,y,k)=\phi(z,y,k)\,\rho_{\l}(z,k)\,,
\eeq
where the coefficients of $\rho_{\l}$ are $\p_U$-invariant holomorphic functions. Then
the same arguments as in \cite{shiota} show that there exists a $\p_U$-invariant series
$f(z,k)$ with holomorphic in $z$ coefficients and formal series $\mu_{\l}^0(k)$ with constant
coefficients such that the equation
\beq\label{v13}
f(z+\l,k)\rho_{\l}(z,k)=f(z,k)\,\mu_{\l}(k)
\eeq
holds. The ambiguity in the choice of $f$ and $\mu$ corresponds to the multiplication by
the exponent of a linear form in $z$ vanishing on $U$, i.e.,
\beq\label{v14}
f'(z,k)=f(z,k)\,e^{(b(k),z)}, \ \ \mu_{\l}'(k)=\mu_{\l}(k)\,e^{(b(k),\l)}, \ \ (b(k),U)=0,
\eeq
where $b(k)=\sum_sb_s k^{-s}$ is a formal series with vector-coefficients that are
orthogonal to $U$. The vector $U$ is in general position with respect to the lattice. Therefore,
the ambiguity can be uniquely fixed by imposing $(d-1)$ normalizing conditions $\mu_{\l_i}(k)=1,\
i>1$ (recall that $\mu_{\l_1}(k)=1$ by construction).

The formal series $f\phi$ is quasi-periodic and its multipliers satisfy
(\ref{v11}). Then,  by that properties it is defined uniquely up to a factor
which is constant in $z$ and $y$.
Therefore, for the unique definition of $\phi_0$ it is enough to fix its evaluation at
$z_0$ and $y=0$. The corollary is proved.

\subsection*{The spectral curve}

The next goal is to show that $\l$-periodic wave solutions of equation (\ref{lax0}), with
$u$ as in (\ref{u0}), are common eigenfunctions of rings of commuting operators.

Note that a simple shift $z\to z+Z$, where $Z\notin \Sigma,$ gives
$\l$-periodic wave solutions with meromorphic coefficients along the affine
subspaces $Z+\mathbb C^d$. Theses $\lambda$-periodic wave solutions are related to each other
by $\p_U$-invariant factor. Therefore choosing, in the neighborhood of any
$Z\notin \Sigma,$ a hyperplane orthogonal to the vector $U$ and
fixing initial data on this hyperplane at $y=0,$ we define the corresponding
series $\phi(z+Z,y,k)$ as a {\it local} meromorphic function of $Z$ and the
{\it global} meromorphic function of $z$.

\begin{lem} Let the assumptions of Theorem~\ref{th1.1} hold. Then there is a unique
pseudo\-differential operator
\beq\label{LL}
\L(Z,\p_x)=\p_x+\sum_{s=1}^{\infty} w_s(Z)\p_x^{-s}
\eeq
such that
\beq\label{kk}
\L(Ux+Vy+Z,\p_x)\,\psi=k\,\psi\,,
\eeq
where $\psi=e^{kx+k^2y} \phi(Ux+Z,y,k)$ is a $\l$-periodic solution of
(\ref{lax0}).
The coefficients $w_s(Z)$ of $\L$  are meromorphic functions on the abelian variety $X$
with poles along  the divisor $\Theta$.
\end{lem}
{\it Proof.}
Let $\psi$ be a $\l$-periodic wave solution. The substitution of (\ref{psi2}) in (\ref{kk})
gives a system of equations
that recursively define $w_s(Z,y)$ as differential polynomials in $\xi_s(Z,y)$.
The coefficients of $\psi$ are local meromorphic functions of $Z$, but
the coefficients of $\L$ are well-defined
{\it global meromorphic functions} of on $\mathbb C^g\setminus\Sigma$, because
different $\l$-periodic wave solutions are related to each other by $\p_U$-invariant
factor, which does not affect $\L$. The singular locus is
of codimension $\geq 2$. Then Hartogs' holomorphic extension theorem implies that
$w_s(Z,y)$ can be extended to a global meromorphic function on $\mathbb C^g$.

The translational invariance of $u$ implies the translational invariance of
the $\l$-periodic wave solutions. Indeed, for any constant $s$ the series
$\phi(Vs+Z,y-s,k)$ and $\phi(Z,y,k)$ correspond to $\l$-periodic solutions
of the same equation. Therefore, they coincide up to a $\p_U$-invariant factor.
This factor does not affect $\L$. Hence, $w_s(Z,y)=w_s(Vy+Z)$.

The $\l$-periodic wave functions corresponding to $Z$ and
$Z+\lambda'$ for any $\lambda'\in \Lambda$
are also related to each other by a $\p_U$-invariant factor:
\beq\label{tr1}
\p_U\left(\phi_1(Z+\l',y,k)\phi^{-1}(Z,y,k)\right)=0.
\eeq
Hence, $w_s$ are periodic with respect to $\Lambda$ and therefore are
meromorphic functions on the abelian variety $X$.
The lemma is proved.

Consider now the differential parts of the pseudo\-differential operators $\L^m$.
Let $\L^m_+$ be the differential operator such that
$\L^m_-=\L^m-\L^m_+=F_m\p^{-1}+O(\p^{-2})$. The leading
coefficient $F_m$ of $\L^m_-$ is the residue of $\L^m$:
\beq\label{res1}
F_m={\rm res}_{\p}\  \L^m.
\eeq
 From the construction of $\L$ it follows that $[\p_y-\p^2_x+u, \L^n]=0$. Hence,
\beq\label{lax}
[\p_y-\p_x^2+u,\L^m_+]=-[\p_y-\p_x^2+u, \L^m_-]=2\p_x F_m
\eeq
(compare with (\ref{kphier4})).
The functions $F_m$ are differential polynomials in the coefficients $w_s$ of $\L$.
Hence, $F_m(Z)$ are meromorphic functions on $X$. Next statement is crucial for
the proof of the existence of commuting differential operators associated with $u$.
\begin{lem} The abelian functions $F_m$ have at most the second order pole on the divisor
$\Theta$.
\end{lem}
{\it Proof.} We need a few more standard constructions from the KP theory.
If $\psi$ is as in Lemma~\ref{lem3.1}, then  there exists a unique pseudo\-differential
operator $\Phi$ such that
\beq\label{S}
\psi=\Phi e^{kx+k^2y},\ \ \Phi=1+\sum_{s=1}^{\infty}\f_s(Ux+Z,y)\p_x^{-s}.
\eeq
The coefficients of $\Phi$ are universal differential polynomials on $\xi_s$.
Therefore, $\f_s(z+Z,y)$ is a global meromorphic function of $z\in C^d$ and
a local meromorphic function of $Z\notin \Sigma$.
Note that $\L=\Phi(\p_x)\, \Phi^{-1}$.

Consider the dual wave function defined by the left action of the operator $\Phi^{-1}$:
$\psi^+=\bigl(e^{-kx-k^2y}\bigr)\Phi^{-1}$.
Recall that the left action of a pseudo\-differential operator is the formal adjoint action
under which the left action of $\p_x$ on a function $f$ is $(f\p_x)=-\p_xf$.
If $\psi$ is a formal wave solution of (\ref{lax0}),
then $\psi^+$ is a solution of the adjoint equation
\beq\label{adj}
(-\p_y-\p_x^2+u)\psi^+=0.
\eeq
The same arguments, as before, prove that if equations (\ref{cm50}) for poles of $u$
hold then $\xi_s^+$ have simple poles at the poles of $u$. Therefore, if $\psi$ is as in
Lemma~\ref{lem6.3}, then the dual wave solution is of the form
$\psi^+=e^{-kx-k^2y}\phi^+(Ux+Z,y,k)$, where
the coefficients $\xi_s^+(z+Z,y)$ of the formal series
\beq\label{psi2+}
\phi^+(z+Z,y,k)=e^{-by}\biggl(1+\sum_{s=1}^{\infty}\xi^+_s(z+Z,y)\, k^{-s}\biggr)
\eeq
are $\l$-periodic meromorphic functions of the variable $z\in \mathbb C^d$ with a
simple pole at the divisor $\Theta^{\,U}(y)$.

The ambiguity in the definition of $\psi$ does not affect the product
\beq\label{J0}
\psi^+\psi=\left(e^{-kx-k^2y}\Phi^{-1}\right)\left(\Phi e^{kx+k^2y}\right).
\eeq
Therefore, although each factor is only a local meromorphic function on
$\mathbb C^g\setminus \Sigma$, the coefficients $J_s$ of the product
\beq\label{J}
\psi^+\psi=\phi^+(Z,y,k)\phi(Z,y,k)=1+\sum_{s=2}^{\infty}J_s(Z,y)k^{-s}.
\eeq
are {\it global meromorphic functions} of $Z$. Moreover, the translational invariance
of $u$ implies that they have the form $J_s(Z,y)=J_s(Z+Vy)$. Each of the factors in
the left-hand side of (\ref{J}) has a simple pole on $\Theta-Vy$. Hence, $J_s(Z)$ is
a meromorphic function on $X$ with a second order pole at $\Theta$. According to Lemma \ref{lem4.3}, we have
$F_n=J_{n+1}$. That completes the proof of the lemma.

Let $\bf {\hat F}$ be a linear space generated by $\{F_m, \ m=0,1,\ldots\}$, where we
set $F_0=1$. It is a subspace of the
$2^g$-dimensional space of the abelian functions that have at most second order pole at
$\Theta$. Therefore, for all but $\hat g=\dim{\bf \hat F}$ positive integers $n$,
there exist constants $c_{i,n}$ such that
\beq\label{f1}
F_n(Z)+\sum_{i=0}^{n-1} c_{i,n}F_i(Z)=0.
\eeq
Let $I$ denote the subset of integers $n$ for which there are no such constants. We call
this subset the gap sequence.
\begin{lem} Let $\L$ be the pseudo\-differential operator corresponding to
a $\l$-periodic wave function $\psi$ constructed above. Then, for the differential operators
\beq\label{a2}
L_n=\L^n_++\sum_{i=0}^{n-1} c_{i,n}\L^{n-i}_+=0, \ n\notin I,
\eeq
the equations
\beq\label{lp}
L_n\,\psi=a_n(k)\,\psi, \ \ \ a_n(k)=k^n+\sum_{s=1}^{\infty}a_{s,n}k^{n-s}
\eeq
where $a_{s,n}$ are constants, hold.
\end{lem}
{\it Proof.} First note that from (\ref{lax}) it follows that
\beq\label{lax3}
[\p_y-\p_x^2+u,L_n]=0.
\eeq
Hence, if $\psi$ is a $\l$-periodic wave solution of (\ref{lax0})
corresponding to $Z\notin \Sigma$, then $L_n\psi$ is also a formal
solution of the same equation. That implies the equation
$L_n\psi=a_n(Z,k)\psi$, where $a$ is $\p_U$-invariant. The ambiguity in the definition of
$\psi$ does not affect $a_n$. Therefore, the coefficients of $a_n$ are well-defined
{\it global} meromorphic functions on $\mathbb C^g\setminus \Sigma$. The $\p_U$-
invariance of $a_n$ implies that $a_n$, as a function of $Z$, is holomorphic outside
of the locus. Hence it has an extension to a holomorphic function on $\mathbb C^g$.
Equations (\ref{tr1}) imply that $a_n$ is periodic with respect to the lattice
$\Lambda$. Hence $a_n$ is $Z$-independent. Note that $a_{s,n}=c_{s,n},\ s\leq n$.
The lemma is proved.

The operator $L_m$ can be regarded as a $Z\notin \Sigma$-parametric family   of
ordinary differential operators $L_m^Z$ whose coefficients have the form
\beq\label{lu}
L_m^Z=\p_x^n+\sum_{i=1}^m u_{i,m}(Ux+Z)\, \p_x^{m-i},\ \ m\notin I.
\eeq
\begin{cor} The operators $L_m^Z$
commute with each other,
\beq\label{com1}
[L_n^Z,L_m^Z]=0, \ Z\notin \Sigma.
\eeq
\end{cor}
 From (\ref{lp}) it follows that $[L_n^Z,L_m^Z]\psi=0$. The commutator is an ordinary
differential operator. Hence, the last equation implies (\ref{com1}).

\begin{lem}\label{lem6.7} Let $\A^Z,\ Z\notin \Sigma,$ be a commutative ring of ordinary differential
operators spanned by the operators $L_n^Z$. Then there is an irreducible algebraic
curve $\G$ of arithmetic genus $\hat g=\dim{\bf \hat F}$ such that $\A^Z$ is isomorphic
to the ring $A(\G,P_0)$ of the meromorphic functions on $\G$ with the only pole at
a smooth point $P_0$. The correspondence $Z\to \A^Z$ defines a holomorphic
imbedding of $X\setminus \Sigma$ into the space of torsion-free rank 1 sheaves $\F$ on $\G$
\beq\label{is}
j\colon X\setminus\Sigma\longmapsto \overline{\rm Pic}(\G).
\eeq
\end{lem}
{\it Proof.} In order to get the statement of the theorem as a direct corollary of Theorem 5.1,
it remains only to show that the ring $\A^Z$ is maximal.
Recall, that a commutative ring $\A$ of linear ordinary differential operators
is called maximal if it is not contained in any bigger commutative ring.
Let us show that for a generic $Z$ the ring $\A^Z$ is maximal. Suppose that it is not.
Then there exits  $\a\in I$, where $I$ is the gap sequence defined above, such
that for each $Z\notin \Sigma$ there exists an operator $L_{\a}^Z$
of order $\a$ which commutes with $L_n^Z, n\notin I$. Therefore, it commutes with $\L$.
A differential operator commuting with $\L$ up to the order $O(1)$
can be represented in the
form $L_{\a}=\sum_{m<\a} c_{i,\a}(Z)\L^i_+$, where $c_{i,\a}(Z)$ are
$\p_1$-invariant functions of $Z$. It commutes with $\L$ if and only if
\beq\label{z3}
F_{\a}(Z)+\sum_{i=0}^{n-1} c_{i,\a}(Z)F_i(Z)=0,\ \ \p_U c_{i,\a}=0.
\eeq
Note the difference between (\ref{f1}) and (\ref{z3}). In the first equation the
coefficients $c_{i,n}$ are constants. The $\l$-periodic wave solution of equation
(\ref{lax0}) is a common eigenfunction of all  commuting operators, i.e.,
$L_{\a}\psi=a_{\a}(Z,k)\psi$, where
$a_{\a}=k^{\a}+\sum_{s=1}^{\infty} a_{s,\a}(Z)k^{\a-s}$ is $\p_1$-invariant.
The same arguments as those used in the proof of equation (\ref{lp}) show that
the eigenvalue $a_{\a}$ is $Z$-independent. We have
$a_{s,\a}=c_{s,\a},\ s\leq \a$. Therefore, the coefficients in (\ref{z3}) are
$Z$-independent. That contradicts the assumption that $\a\notin I$.
The lemma is proved.

Our next goal is to prove finally the global existence of the wave function.
\begin{lem} Let the assumptions of the Theorem~\ref{th1.2} hold. Then there exists a common
eigenfunction
of the corresponding commuting operators $L_n^Z$ of the form
$\psi=e^{kx}\phi(Ux+Z,k)$ such that
the coefficients of the formal series
\beq\label{psi6}
\phi(Z,k)=1+\sum_{s=1}^{\infty}\xi_s(Z)\, k^{-s}
\eeq
are global meromorphic functions with a simple pole at $\Theta$.
\end{lem}
{\it Proof.} It is instructive to consider first the case when the spectral curve $\G$
of the rings $\A^Z$ is smooth. Then, as shown in (\cite{kr1, kr2}),
the corresponding common
eigenfunction of the commuting differential operators (the Baker-Akhiezer function),
normalized by the
condition $\psi_0|_{x=0}=1$, is of the form
(\cite{kr1, kr2})
\beq\label{ba}
\hat \psi_{0}={\hat\theta (\hat A(P)+\hat Ux+\hat Z)\,\hat\theta (\hat Z)\over
\hat\theta(\hat Ux+\hat Z)\,\hat\theta(\hat A(P)+\hat Z)}\,
e^{x\,\Omega(P)}.
\eeq
(compare with (\ref{BA1}). Here $\hat \theta (\hat Z)$ is the Riemann theta-function constructed
with the help of the
matrix of $b$-periods of normalized holomorphic differentials on $\G$; $\hat A\colon \G\to J(\G)$
is the Abel-Jacobi map; $\Omega$ is the abelian integral corresponding to the second kind
meromorphic differential $d\Omega$ with the only pole of the form $dk$ at the marked point
$P_0$ and $2\pi i \hat U$ is the vector of its $b$-periods.

\noindent{\bf Remark.}
Let us emphasize, that the formula (\ref{ba}) is not the result of solution of some
differential equations.
It is a direct corollary of analytic properties of the Baker-Akhiezer function
$\hat \psi_0(x,P)$ on the spectral curve.

The last factors in the numerator and the denominator of (\ref{ba}) are $x$-independent.
Therefore, the function
\beq\label{ba1}
\hat \psi_{BA}={\hat\theta (\hat A(P)+\hat Ux+\hat Z)\over
\hat\theta(\hat Ux+\hat Z)}\,
e^{x\,\Omega(P)}
\eeq
is also a common eigenfunction of the commuting operators.

In the neighborhood of $P_0$ the function $\hat \psi_{BA}$ has the form
\beq\label{sh10}
\hat \psi_{BA}=e^{kx}
\biggl(1+\sum_{s=1}^{\infty}{\tau_s (\hat Z+\hat Ux)\over \hat \theta(\hat U x+\hat Z)}\,k^{-s}
\biggr), \ \ k=\Omega,
\eeq
where $\tau_s(\hat Z)$ are global holomorphic functions.

According to Lemma~\ref{lem6.7}, we have a holomorphic imbedding $\hat Z=j(Z)$
of $X\setminus\Sigma$ into $J(\G)$.
Consider the formal series $\psi=j^*\hat \psi_{BA}$.
It is globally well-defined out of $\Sigma$.
If $Z\notin \Theta$, then $j(Z)\notin \hat \Theta$
(which is the divisor on which the condition (\ref{sheaf}) is violated).
Hence, the coefficients of $\psi$ are regular out of $\Theta$. The singular locus
is at least of codimension 2. Hence, using once again Hartogs' arguments we can
extend $\psi$ on $X$.

If the spectral curve is singular, we can proceed along the same lines
using the generalization of (\ref{ba1})
given by the theory of Sato $\tau$-function
(\cite{wilson}). Namely, a set of algebraic-geometrical data (\ref{corr}) defines
the point of the Sato Grassmannian, and therefore, the corresponding $\tau$-function: $\tau(t;\F)$.
It is a holomorphic function of the variables $t=(t_1,t_2,\ldots)$,
and is a section of a holomorphic line bundle on $\overline{\rm Pic}(\G)$.

The variable $x$ is identified with the first time of the KP-hierarchy, $x=t_1$.
Therefore, the formula for the Baker-Akhiezer function corresponding to a point of the Grassmannian
(\cite{wilson}) implies that the function $\hat \psi_{BA}$ given by the formula
\beq\label{baw}
\hat \psi_{BA}={\tau(x-k,-{1\over 2}k^2,-{1\over 3}k^3,\ldots; \F) \over
\tau(x,0,0,\ldots; \F)}e^{kx}
\eeq
is a common eigenfunction of the commuting operators defined by $\F$.
The rest of the arguments proving the lemma are
the same, as in the smooth case.
\begin{lem} The linear space $\bf {\hat F}$ generated by the abelian functions $\{F_0=1,
F_m=\res_\p \L^m\},$ is
a subspace of the space $\bf H$ generated by $F_0$ and by the abelian functions
$H_i=\p_U\p_{z_i}\ln \theta(Z)$.
\end{lem}
{\it Proof.} Recall that the functions $F_n$ are abelian
functions with at most second order pole on $\Theta$. Hence, a~priori
$\hat g=\dim{\bf \hat F}\leq 2^g.$ In order to prove the statement of the lemma
it is enough to show that $F_n=\p_U Q_n$, where $Q_n$ is a
meromorphic function with a pole along $\Theta$. Indeed, if $Q_n$ exists,
then, for any vector $\l$ in the period lattice, we have $Q_n(Z+\l)=Q_n(Z)+c_{n,\l}$.
There is no abelian function with a simple pole on $\Theta$. Hence, there exists
a constant $q_n$ and two $g$-dimensional vectors $l_n,l_n'$, such that
$Q_n=q_n+(l_n,Z)+(l_n',h(Z))$, where $h(Z)$ is a vector with the coordinates
$h_i=\p_{z_i}\ln \theta$. Therefore, $F_n=(l_n,U)+(l_n',H(Z))$.

Let $\psi(x,Z,k)$ be the formal Baker-Akhiezer function defined in the previous lemma.
Then the coefficients $\varphi_s(Z)$ of the corresponding wave operator $\Phi$ (\ref{S})
are global meromorphic functions with poles on $\Theta$.

The left and  right action of pseudo\-differential operators are formally adjoint,
i.e., for any two operators the equality $\left(e^{-kx}\D_1\right)\left(\D_2e^{kx}\right)=
e^{-kx}\left(\D_1\D_2e^{kx}\right)+\p_x\left(e^{-kx}\left(\D_3e^{kx}\right)\right)$
holds. Here $\D_3$ is a pseudo\-differential operator whose coefficients are differential
polynomials in the coefficients of $\D_1$ and $\D_2$. Therefore, from (\ref{J0})
it follows that
\beq\label{z8}
\psi^+\psi=1+\sum_{s=2}^{\infty}F_{s-1}k^{-s}=
1+\p_x\biggl(\sum_{s=2}^{\infty}Q_sk^{-s}\biggr).
\eeq
The coefficients of the series $Q$ are differential polynomials in the
coefficients $\varphi_s$ of the wave operator. Therefore, they are global meromorphic
functions of $Z$ with poles on $\Theta$. Lemma is proved.

The construction of multivariable \BA functions presented in Section 2 for smooth curves
is a manifestation of general statement valid for singular spectral curves:
flows of the KP hierarchy define deformations of the commutative
rings $\A$ of ordinary linear differential operators. The spectral curve
is invariant under these flows. For a given spectral curve $\G$ the orbits of
the KP hierarchy are isomorphic to the generalized Jacobian $J(\G)={\rm Pic}^0 (\G)$,
which is the equivalence classes of zero degree divisors on the spectral curve
(see details in \cite{shiota,kr1,kr2,wilson}).

As shown in Section 4, the evolution of the potential $u$ is described by equation (\ref{kphier})
The first two times of the hierarchy are identified with the variables $t_1=x, t_2=y$.
Equations (\ref{kphier}) identify the space $\hat {\bf F}_1$ generated by the functions
$\p_U F_n$ with the tangent space of the KP orbit at  $\A^Z$.
Then, from Lemma 6.9 it follows that this tangent space is a subspace of the tangent space
of the abelian variety $X$. Hence, for any $Z\notin \Sigma$, the orbit of the KP flows
of the ring $\A^Z$ is in $X$, i.e., it defines an holomorphic imbedding:
\beq\label{imb}
i_Z\colon J(\G)\longmapsto X.
\eeq
 From (\ref{imb}) it follows that $J(\G)$ is {\it compact}.

The generalized Jacobian of an algebraic curve is compact
if and only if the curve is {\it smooth} (\cite{mdl}). On a smooth algebraic curve
a torsion-free rank 1 sheaf is a line bundle, i.e., $\overline {\rm Pic} (\G)=J(\G)$.
Then (\ref{is}) implies that $i_Z$ is an isomorphism. Note that
for the Jacobians of smooth algebraic curves the bad locus $\Sigma$ is empty
(\cite{shiota}), i.e., the imbedding $j$ in (\ref{is}) is defined everywhere
on $X$ and is inverse to $i_Z$. Theorem~\ref{th1.1} is proved.

\section{Characterization of the Prym varieties}

To begin with let us recall the definition of Prym varieties.
An involution $\s\colon \G\longrightarrow \G$ of a smooth algebraic curve $\G$
induces an involution $\s^*\colon  J(\G)\longrightarrow J(\G)$ of the Jacobian.
The kernel of the map $1+\s^*$ on $J(\G)$ is the sum of a lower-dimensional
abelian variety, called the Prym variety (the connected component of
zero in the kernel), and a finite group. The Prym variety
naturally has a polarization induced by the principal polarization
on $J(\G)$. However, this polarization is not principal, and
the Prym variety admits a natural principal polarization if and only
if $\s$ has at most two fixed points on $\G$ --- this is the case we will concentrate on.

 From the point of view of integrable systems, attempts to prove the analog of Novikov's conjecture for the case of Prym varieties of algebraic curves with two smooth fixed points of involution were made in \cite{taim2,shiotaPrym,bd}. In \cite{taim2} it was shown that Novikov-Veselov (NV) equation provides solution of the characterization problem up to possible existence of additional irreducible components.
In \cite{shiotaPrym,bd} the characterizations of the Prym varieties in terms of BKP and NV equations were proved only under certain additional assumptions. Moreover, in \cite{bd} an example of a ppav that is not a Prym but for which the theta function gives a solution to the BKP equation was constructed. Thus for more than 15 years it was
widely accepted that Prym varieties can not be characterized with the help if integrable systems.

In \cite{kr-prym} the first author proved that Prym varieties of algebraic curves with two smooth fixed points of involution are characterized among all ppavs by the property of their theta functions providing explicit formulas {\it for solutions of the integrable $2D$ Schr\"odinger equation}, which is one of the auxiliary linear problems for the Novikov-Veselov equation.

Prym  varieties possess generalizations of some properties of Jacobians. In \cite{bd} Beauville and Debarre, and in \cite{fay2} Fay showed that the Kummer images of Prym varieties admit a 4-dimensional family of quadrisecant
planes (as opposed to a 4-dimensional family of trisecant lines for Jacobians). Similarly to the case of Jacobians, it was then shown by Debarre in \cite{deb}
that the existence of a one-dimensional family of quadrisecants characterizes Prym
varieties among all ppavs. However, Beauville and Debarre in \cite{bd} constructed a
ppav that is {\it not} a Prym but such that its Kummer image {\it has} a quadrisecant
plane. Thus no analog of the trisecant conjecture for Prym varieties was conjectured,
and the question of characterizing Prym varieties by a finite amount of geometric data
(i.e., by polynomial equations for theta functions at a finite number of points)
remained completely open.

In \cite{kr-quad} S. Grushevsky and the first author proved that Prym varieties
of unramified covers are characterized among all ppavs by the property of their Kummer images admitting a {\it symmetric pair of quadrisecant 2-planes}.
That there exists such a symmetric pair of quadrisecant planes for the Kummer image of
a Prym variety can be deduced from the description of the 4-dimensional family of
quadrisecants, using the natural involution on the Abel-Prym curve. However, the
statement that a symmetric pair of quadrisecants in fact characterizes Pryms seems
completely unexpected.

The geometric characterization of Prym varieties follows from a characterization of
Prym varieties among all ppavs by some theta-functional equations,
which by using Riemann's bilinear addition
theorem can be shown to be equivalent  to the existence of a symmetric pair of
quadrisecant planes. In order to obtain such a characterization of Prym varieties in \cite{kr-quad}
a {\it new} hierarchy of difference equations, starting from a discrete version of the Schr\"odinger equation was introduced, developed, and studied . The hierarchy constructed can be thought of as a
discrete analog of the Novikov-Veselov hierarchy.

\begin{theo}[Main theorem]\label{theoremPrym}
An indecomposable principally polarized abelian variety $(X,\theta)\in\A_g$
lies in the closure of the locus $\P_g$ of Prym varieties of
unramified double covers if and only if there exist vectors
$A,U,V, W\in\mathbb C^g$ representing distinct points in $X$, none of them points of order two, and constants
$c_1,c_2,c_3,w_1,w_2,w_3\in\mathbb C$ such that one of the following
equivalent conditions holds:\smallskip

$(A)$  The difference $2D$ Schr\"odinger equation
\begin{equation}\label{laxdd-prym}
  \psi_{n+1,m+1}-u_{n,m}(\psi_{n+1,m}-\psi_{n,m+1})-\psi_{n,m}=0,
\end{equation}
with
\begin{equation}\label{ud-prym}
  u_{n,m}:=C_{nm} {\theta((n+1)U+mV+\nu
  W+Z)\,\theta(nU+(m+1)V+\nu W+Z)\over
  \theta((n+1)U+(m+1)V+\overline\nu W+Z)\,\theta(nU+mV+\overline\nu W+Z)}
\end{equation}
and
\begin{equation}\label{pd-prym}
  \psi_{n,m}:={\theta(A+nU+mV+\nu_{nm} W+Z)\over
  \theta(nU+mV+\overline\nu_{nm} W+Z)}\, w_1^nw_2^mw_3^{\nu_{nm}}
  \left(c_1^mc_2^n\right)^{1-2\nu_{nm}},
\end{equation}
is satisfied for all $Z\in X$, where
\begin{equation}\label{ud11}
  \nu:=\nu_{nm}:={1+(-1)^{n+m+1}\over2},\ \ \,
  \overline\nu:=1-\nu, \ \ \,
  C_{nm}:=c_3\left(c_2^{2n+1}c_1^{2m+1}\right)^{1-2\nu_{nm}}\!.
\end{equation}
\smallskip

$(B)$ The following identity holds:
\begin{multline*} 
w_1w_2(c_1c_2)^{\pm 1}\wt K\left({A+U+V\mp W\over
2}\right)-w_1c_3(w_3c_1)^{\pm1}\wt K\left({A+U-V\pm W\over 2}\right)\\
+w_2c_3(w_3c_2)^{\pm 1}\wt K\left({A+V-U\pm W \over 2}\right)-
\wt K\left({A-U-V\mp W\over 2}\right) =0\,,
\end{multline*}
where $\wt K\colon \bC^g \ni z\mapsto
\bigl(\Theta[\e,0](z)\bigr)\in\bC^{2^g}$
is a lifting of the Kummer map (\ref{kum}) to the universal covering of $X$.
\smallskip

$(C)$ The two equations (one for the top choice of signs everywhere,
and one for the bottom)
\begin{multline}\label{cm7d-prym}
c_1^{\mp 2}c_3^2\ \theta(Z+U-V)\,\theta(Z-U\pm W)\,\theta(Z+V\pm W)\\
+c_2^{\mp 2}c_3^2\ \theta(Z-U+V)\,\theta(Z+U\pm W)\,\theta(Z-V\pm W)\\
=c_1^{\mp 2}c_2^{\mp 2}\,\theta(Z-U-V)\,\theta(Z+U\pm W)\,\theta(Z+V\pm W)\\
+\theta(Z+U+V)\,\theta(Z-U\pm W)\,\theta(Z-V\pm W)
\end{multline}
are valid on the theta divisor $\{Z\in X: \theta(Z)=0\}$.
\end{theo}
A purely geometric restatement of part $(B)$ of this result is as
follows.
\begin{cor}[Geometric characterization of Pryms]
A ppav $(X,\theta)\in\A_g$ lies in the closure of the locus of Prym varieties of unramified (\'etale) double covers if and only there exist four distinct points $p_1,p_2,p_3,p_4\in
X$, none of them points of order two, such that the following two quadruples of points on the Kummer variety of $X$:
$$
\{ K(p_1+\e_2p_2+\e_3p_3+\e_4p_4)\mid \e_i\in\{\pm1\},\ \e_2\e_3\e_4=+1\}
$$
and
$$
\{ K(p_1+\e_2p_2+\e_3p_3+\e_4p_4)\mid \e_i\in\{\pm1\},\ \e_2\e_3\e_4=-1\}
$$
are linearly dependent.

Equivalently, this can be stated as saying that $(X,\theta)$ lies in the closure of the Prym if and only if there exists a pair of symmetric (under the $z\mapsto 2p_1-z$ involution) quadrisecants of $K(X)$.
\end{cor}

At first glance the structure of the proof is the same as above. It begins with
a construction of a wave solution of the discrete
analog of $2D$ Schr\"odinger equation (\ref{laxdd-prym}). But in fact, the hierarchy considered involves essentially a pair of functions and is thus essentially a matrix hierarchy, unlike the scalar hierarchy arising for the trisecant case. The argument is very delicate, and involves using the pair of quadrisecant conditions to recursively construct a pair of auxiliary solutions (essentially corresponding to the two components of the kernel, only one of which is the Prym). We refer the reader to \cite{kr-quad}) for details.

Our goal for this section is to elaborate on the ``only if'' part of the statement of the theorem, because
as a byproduct it gives new identities for theta-function which are poorly understood an seems require additional attention.

\subsection*{Four point Baker-Akhiezer function} Four-point Baker-Akhiezer function depends on three discrete
parameters and, as was mentioned in Section 2 gives solution to the BDHE equation. For various choice of
two linear combination of these variables one obtain various linear equation. In \cite{krdd} (see
details in \cite{kwz}) it was shown that the following choice of the ``discrete times'' gives a
a construction of algebraic-geometric 2D difference Schr\"odinger operators.

Let $\G$ be a smooth algebraic curve of genus $\hat g$. Fix four
points $P_1^{\, \pm}, P_2^{\, \pm}\in \G$, and let $\hat
D=\g_1+\cdots+\g_{\wh g}$ be a generic effective divisor on $\G$ of
degree $\wh g$.  By the Riemann-Roch theorem one computes $h^0(\hat
D+n(P_1^+-P_1^-)+m(P_2^+-P_2^-))=1$, for any $n,m\in\mathbb Z$, and for
$\hat D$ generic. We denote by $\wh\psi_{n,m}(P),\ P\in \G$ the
unique section of this bundle. This means that $\wh\psi_{n,m}$ is the unique up to a constant factor meromorphic function such that (away from the marked points $P_i^{\pm}$) it has poles only at $\g_s$, of
multiplicity not greater than the multiplicity of $\g_s$ in $\wh D$, while at the points $P_1^+, P_2^+$ (resp.\ $P_1^-,P_2^-$) the function $\wh\psi_{n,m}$ has poles (resp.\ zeros) of orders $n$ and
$m$.

If we fix local coordinates $k^{-1}$ in the neighborhoods of
marked points (it is customary in the subject to think of marked
points as punctures, and thus it is common to use coordinates such
that $k$ at the marked point is infinite rather than zero), then
the Laurent series for $\psi_{n,m}(P)$, for $P\in\G$ near a marked
point, has the form
\begin{align}
\wh\psi_{n,m}&=k^{\pm n}\biggl(\sum_{s=0}^{\infty}\xi_s^{\,
\pm}(n,m)k^{-s}\biggr), \ \ k=k(P), \
P\to P_1^{\, \pm}, \label{2}\\
\wh\psi_{n,m}&=k^{\pm m}\biggl(\sum_{s=0}^{\infty}\chi_s^{\,
\pm}(n,m)k^{-s}\biggr), \ \ k=k(P), \ P\to P_2^{\, \pm}. \label{3}
\end{align}
As it was shown in Section 2 the function $\psi_{n,m}$ can be expressed
as follows:
\begin{equation}\label{4}
\wh\psi_{n,m}(P)=r_{nm}{\wh\theta(\wh A(P)+n\wh U+m\wh V+\wh Z)
\over \wh \theta(\wh A(P)+\wh Z)\,}\
e^{n\wh\Omega_1(P)+m\wh\Omega_2(P)},
\end{equation}
where  for $i=1,2$ the differential
$d\wh\Omega^i\in H^0(K_\G+P_i^++P_i^-)$ is of the
third kind, normalized to have residues $\mp 1$ at $P_i^{\, \pm}$ and with zero integrals over all the $a$-cycles,
and $\wh\Omega^i$ is the corresponding abelian integral;
we have the following expression
$r_{nm}$ is some constant, $\wh U=\wh A(P_1^-)-\wh A(P_1^+), \
\ \wh V=\wh A(P_2^-)-\wh A(P_2^+),$ and
\begin{equation}\label{hatz}
\wh Z=-\sum_s \wh A(\g_s)+\wh\kappa,
\end{equation}
where $\wh\kappa$ is the vector of Riemann constants.

\paragraph*{Change of notation} We use here notation $\wh \theta$ for the Riemann theta-function of $\G$, for later use of $\theta$ for the Prym theta function.

\begin{theo}[\cite{krdd}]
The Baker-Akhiezer function $\wh \psi_{n,m}$ given by formula (\ref{4})
satisfies the following difference equation
\begin{equation}\label{eqn}
\wh\psi_{n+1,m+1}-a_{n,m}\wh\psi_{n+1,m}-b_{n,m}\wh\psi_{n,m+1}
+c_{n,m}\wh\psi_{n,m}=0,
\end{equation}
\end{theo}

\subsection*{Setup for the Prym construction}

We now assume that the curve $\G$ is an algebraic curve endowed with
an involution $\s$ without fixed points; then $\G$ is a unramified
double cover $\G\longrightarrow \G_0$, where $\G_0=\G/\s$. If $\G$ is of
genus $\wh g=2g+1$, then by Riemann-Hurwitz the genus of $\G_0$ is
$g+1$. From now on we assume that $g>0$ and thus $\wh g>1$. On $\G$ one can choose a basis of cycles $a_i, b_i$ with the
canonical matrix of intersections $a_i\cdot a_j=b_i\cdot b_j=0,\
a_i\cdot b_j=\delta_{ij},\ \ 0\leq i,j\leq 2g,$ such that under the
involution $\s$ we have $\s(a_0)=a_0,\ \s(b_0)=b_0,\
\s(a_j)=a_{g+j},\ \s(b_j)=b_{g+j}, 1\leq j\leq g$. If $d\omega_i$
are normalized holomorphic differentials on $\G$ dual to this choice
of $a$-cycles, then the differentials
$du_j=d\omega_j-d\omega_{g+j}$, for $j=1\ldots g$ are odd, i.e.,
satisfy $\s^*(du_k)=-du_k$, and we call them the normalized
holomorphic Prym differentials. The matrix of their $b$-periods
\begin{equation}\label{pi}
  \Pi_{kj}=\oint_{b_k}du_j,\ \ 1\leq k,j\leq g\,,
\end{equation}
is symmetric, has positive definite imaginary part, and defines the
Prym variety $$\P(\G):=\mathbb C^g/(\mathbb Z^g+\Pi\mathbb Z^g)$$ and the
corresponding Prym theta function $$\theta(z):=\theta(z,\Pi),$$ for
$z\in\mathbb C^g$. We assume that the marked points $P_1^{\, \pm},
P_2^{\, \pm}$ on $\G$ are permuted by the involution, i.e.,
$P_i^+=\s(P_i^-)$. For further use let us fix in addition a third
pair of points $P_3^{\pm},$ such that also $P_3^-=\s(P_3^+)$.

The Abel-Jacobi map $\G\hookrightarrow J(\G)$ induces the Abel-Prym map
$A\colon \G\longrightarrow\P(\G)$ (this is the composition of the Abel-Jacobi map $\widehat
A\colon \g\hookrightarrow J(\G)$ with the projection $J(\G)\to\P(\G)$). There is a choice of
the base point involved in defining the Abel-Jacobi map, and thus in the Abel-Prym
map; let us choose this base point (such a choice is unique up to a point of order two
in $\P(\G)$) in such a way that
\begin{equation}\label{prab}
  A(P)=-A(\s(P)).
\end{equation}

\subsection*{Admissible divisors} An effective divisor on $\G$
of degree $\hat g-1=2g$, $D=\g_1+\ldots\g_{2g}$, is called {\it
admissible} if it satisfies
\begin{equation}\label{const}
  [D]+[\s(D)]=K_\G\in J(\G)
\end{equation}
(where $K_\G$ is the canonical class of $\G$), and if moreover
$H^0(D+\s(D))$ is generated by an even holomorphic differential
$d\Omega$, i.e., that
\begin{equation}\label{const0}
  d\Omega(\g_s)=d\Omega(\s(\g_s))=0,\ \  d\Omega=\s^*(d\Omega).
\end{equation}

Algebraically, what we are saying is the following. The divisors $D$
satisfying (\ref{const}) are the preimage of the point $K_\G$ under
the map $1+\s$, and thus are a translate of the subgroup
$Ker(1+\s)\subset J(\G)$ by some vector. As shown by Mumford
\cite{mum2}, this kernel has two components --- one of them being
the Prym, and the other being the translate of the Prym variety by
the point of order two corresponding to the cover $\Gamma\to
\Gamma_0$ as an element in $\pi_1(\Gamma_0)$. The existence of an
even differential as above picks out one of the two components, and
the other one is obtained by adding $A-\sigma(A)$ to the divisor of
such a differential, for some $A$.
statement.

\begin{prop}\label{theta.adm}
For a generic vector $Z$ the zero-divisor $D$ of the function
$\theta(A(P)+Z)$ on $\G$ is of degree $2g$ and satisfies the
constraints (\ref{const}) and (\ref{const0}), i.e., is admissible.
\end{prop}
{\bf Remark.} S. Grushevsky and the first author had been unable to find a complete proof of
precisely this statement in the literature. However, both Elham Izadi and Roy Smith have
independently supplied them with simple proofs of this result, based
on Mumford's description and results on Prym varieties.
As pointed out by a referee, this result can also be easily obtained by applying
Fay's proposition 4.1 in \cite{fay}. In \cite{kr-quad} independent analytic proof
was proposed which also can be seen  analytic proof of some of Mumford's
results.

Note that the function $\theta(A(P)+Z)$ is multi-valued on $\G$, but
its zero-divisor is well-defined. The arguments identical to that in
the standard proof of the inversion formula (\ref{hatz}) show that
the zero divisor $D(Z):=\theta(A(P)+Z)$ is of degree $\hat g-1=2g$.

\begin{lem}\label{dOi}
For any pair of points $P_j^\pm$ conjugate under the involution $\sigma$ there exists
a unique differential $d\Omega_j$ of the third kind (i.e., a dipole differential with
simple poles at these points and holomorphic elsewhere), such that it has residues
$\mp 1$ at these points, is odd under $\s$, i.e.,  satisfies
$d\Omega_j=-\s^*(d\Omega_j)$, and such that all of its $a$-periods are integral
multiples of $\pi i$, i.e., such a differential $d\Omega_i$ exists for a unique set of
numbers $l_0,\ldots,l_{g}\in\mathbb Z$ satisfying
\begin{equation}\label{a-cycle}
\oint_{a_k}d\Omega_j=\pi i \,l_k, \ \ k=0,\ldots,g.
\end{equation}
\end{lem}
Indeed, by Riemann's bilinear relations there exists a unique differential
$d\Omega$ of the third kind with residues as required, and satisfying
$\oint_{a_k}d\Omega=0$ for all $k$. Note, however, that then
$\oint_{a_k}\sigma^*(d\Omega)$ is not necessarily zero, as the image $\s(a_k)$ of the
loop $a_k$, while homologous to $a_{g+k}$ on $\tilde{\G}$, is not necessarily homologic
to $a_{g+k}$ (resp.\ to $a_0$ for $\sigma(a_0)$) on $\tilde{\G}\setminus\lbrace
P_j^\pm\rbrace$. Thus each integral $\oint_{a_k}\sigma^*(d\Omega)$, being equal to
$2\pi i$ times the winding number of $\sigma(a_k)$ around $P_j^+$ minus that around
$P_j^-$, is equal to $2\pi i l_k$ for some $l_k\in\mathbb Z$. We now subtract from
$d\Omega$ the linear combination $\pi i\left(l_0d\omega_0+\sum_{k=1}^g
l_k(d\omega_k+d\omega_{g+k})\right)$ of even abelian differentials to get the desired
$d\Omega_j$.

\begin{theo} \cite{kr-quad}
For a generic $D=D(Z)$ and for each set of integers $(n,m,r)$
such that
\begin{equation}\label{ev}
n+m+r=0 \  {\rm mod}\  2
\end{equation}
the space
$$
  H^0(D+n(P_1^+-P_1^-)+m(P_2^+-P_2^-)+r(P_3^+-P_3^-))
$$
is one-dimensional. A basis element of this space is given by
\begin{equation}\label{psiprym}
\psi_{n,m,r}(P):=h_{n,m,r}{\theta(A(P)+nU+mV+rW+Z)\over\theta( A(P)+Z)}\
e^{n\Omega_1(P)+m\Omega_2(P)+r\Omega_3(P)},
\end{equation}
where $\Omega_j$ is the abelian integral corresponding to the differential $d\Omega_j$
defined by lemma~\ref{dOi}, and $U$, $V$, $W$ are the vectors of $b$-periods of these
differentials, i.e.,
\begin{equation}\label{per}
2\pi i U_k=\oint_{b_k}d\Omega_1,\
2\pi i V_k=\oint_{b_k}d\Omega_2,\
2\pi i W_k=\oint_{b_k}d\Omega_3.
\end{equation}
\end{theo}
The proof is identically the same as the proof of (\ref{2.101}).
It is easy to check that the right-hand side of (\ref{psiprym})
is a single valued function on $\G$ having all the desired
properties, and thus it gives a section of the desired bundle. Note
that the constraint (\ref{ev}) is required due to (\ref{a-cycle}),
and the uniqueness of $\psi$ up to a constant factor, i.e., the
one-dimensionality of the $H^0$ above, is a direct corollary of the
Riemann-Roch theorem.

Note that bilinear Riemann identities imply
\begin{equation}\label{periods}
  2U=A(P_1^-)-A(P_1^+), \  \ 2V=A(P_2^-)-A(P_2^+), \ \ 2W=A(P_3^-)-A(P_3^+).
\end{equation}

Let us compare the definition of $\wh \psi_{n,m}$ defined for any
curve $\G$, with that of $\psi_{n,m,r}$, which is only defined for a
curve with an involution satisfying a number of conditions. To make
such a comparison, consider the divisor $\wh D=D+P_3^+$ of degree
$\hat g=2g+1$, and let $\wh\psi_{n,m}$ be the corresponding
Baker-Akhiezer function.
\begin{cor}
For the Baker-Akhiezer function $\wh \psi_{nm}$ corresponding to the
divisor $\wh D=D+P_3^+$ we have
\begin{equation}\label{psi-psi}
  \wh \psi_{nm}=\psi_{n,m,\nu}
\end{equation}
where $\nu=\nu_{nm}$ is defined in (\ref{ud11}), i.e., is 0 or 1 so that
$n+m+\nu$ is even.
\end{cor}
\begin{cor}
If $n+m$ is even, then by formulae (\ref{4}), (\ref{psiprym})
\begin{multline}\label{comp}
 {\wh\theta(\wh A(P)+n\wh U+m\wh V+\wh Z)
  \,\wh \theta(\wh A(P_0)+\wh Z)\over
  \wh \theta(\wh A(P)+\wh Z)\,
  \wh\theta(\wh A(P_0)+n\wh U+m\wh V+\wh Z)}= \\
  {\theta(A(P)+nU+mV+ Z)\, \theta( A(P_0)+Z)\over
  \theta( A(P)+Z)\, \theta(A(P_0)+nU+mV+ Z)}e^{nr_1+mr_2},
\end{multline}
where $r_i=\int_{P_0}^P(d\wh \Omega_i-d\Omega_i)$, and we recall
that $\wh Z=\wh A(\wh D)+\wh \kappa$, and $Z$ is its image.
\end{cor}
{\bf Remark.} This equality, valid for any pair of points $P,P_0$ is a nontrivial
identity between theta functions. The first author's attempts to derive it
directly from the Schottky-Jung relations have failed so far.

\paragraph*{Notation}
For brevity throughout the rest of the paper we use the notation:
$\psi_{n,m}:=\psi_{n,m,\nu_{nm}}$.

\begin{lem}\cite{kr-quad}
The Baker-Akhiezer function $\psi_{n,m}$ given by
\begin{equation}\label{psipr}
\psi_{n,m}={\theta(A(P)+Un+Vm+\nu_{nm} W+Z)\over
  \theta(Un+Vm+\overline\nu_{nm} W + Z)\, \theta( A(P)+Z)}
  \cdot{e^{n\Omega_1(P)+m\Omega_2(P)+\nu_{nm}\Omega_3(P)}\over
 e^{(2\nu_{nm}-1)(n\Omega_1(P_3^+)+m\Omega_2(P_3^+))}},
\end{equation}
where $\overline\nu_{nm}=1-\nu_{nm}$ as in (\ref{ud11}),
satisfies the equation (\ref{laxdd-prym}), i.e.,
$$
  \psi_{n+1,m+1}-u_{n,m}(\psi_{n+1,m}-\psi_{n,m+1})-\psi_{n,m}=0,
$$
with $u_{n,m}$ as in (\ref{ud-prym}), (\ref{ud11}), where
\begin{equation}\label{d2}
  c_1=e^{\Omega_2(P_3^+)},\ \  c_2=e^{\Omega_1(P_3^+)},\ \ c_3=e^{\Omega_1(P_2^+)}
\end{equation}
\end{lem}
Note that the first and the last factors in the denominator of
(\ref{psipr}) correspond to a special choice of the normalization
constants $h_{n,m,\nu}$ in (\ref{psiprym}):
\begin{equation}\label{nordec}
\begin{split}
  \psi_{nm}(P_3^-)&=(\theta(Z+W))^{-1}, \ \ \nu_{nm}=0,\\
  \psi_{nm}e^{-\Omega_3}|_{\,P=P_3^+}&=(\theta(Z-W))^{-1} ,\ \ \nu_{nm}=1.
\end{split}
\end{equation}
This normalization implies that for even $n+m$ the difference
$(\psi_{n+1,m+1}-\psi_{n,m})$ equals zero at $P_3^-$. At the same
time as a corollary of the normalization we get that
$(\psi_{n+1,m}-\psi_{n,m+1})$ has no pole at $P_3^+$. Hence, these
two differences have the same analytic properties on $\G$ and thus
are proportional to each other (the relevant $H^0$ is
one-dimensional by Riemann-Roch). The coefficient of proportionality
$u_{nm}$ can be found by comparing the singularities of the two
functions at $P_1^+$.

The second factor in the denominator of the formula (\ref{psipr})
does not affect equation (\ref{laxdd-prym}). Hence, the lemma proves the
``only if'' part of the statement $(A)$ of the main theorem for the case
of smooth curves. It remains valid under degenerations to singular
curves which are smooth outside of fixed points $Q_k$ which are
simple double points, i.e., to the curves of type $\{\G,\s,Q_k\}$.

\noindent{\bf Remark.}
Equation (\ref{laxdd-prym}) as a special reduction of (\ref{eqn}) was
introduced in \cite{grin}. It was shown that equation (\ref{eqn})
implies a five-term equation
\begin{equation}\label{51}
  \psi_{n+1,m+1}-\tilde a_{nm}\psi_{n+1,m-1}-\tilde b_{n,m}\psi_{n-1,m+1}+
  \tilde c_{nm} \psi_{n-1,m-1}
  =\tilde d_{n,m}\psi_{n,m}
\end{equation}
if and only if it is of the form (\ref{laxdd-prym}). A reduction of the
algebro-geometric construction proposed in \cite{krdd} in the case of
algebraic curves with involution having two fixed points was found.
It was shown that the corresponding Baker-Akhiezer functions do
satisfy an equation of the form (\ref{laxdd-prym}). Explicit formulae for
the coefficients of the equations in terms of Riemann
theta-functions were obtained. The fact that the Baker-Akhiezer
functions and the coefficients of the equations can be expressed in
terms of Prym theta-functions was first obtained in \cite{kr-quad}.

The statement that $\psi_{n,m}$ satisfy (\ref{51}) can be proved
directly. Indeed all the functions involved in the equation are in
$$
  H^0(D+(n+1)P_1^+-(n-1)P_1^-+(m+1)P_2^+-(m-1)P_2^-+\nu(P_3^+-P_3^-))
$$
By the Riemann-Roch theorem the dimension of the latter space is $4$. Hence,
any five elements of this space are linearly dependent, and it remains to
find the coefficients of (\ref{51}) by a
comparison of singular terms at the points $P_1^{\pm}, P_2^{\pm}$.

\begin{theo}\cite{kr-quad}
For any four points $A,U,V,W$ on the image
$\G\hookrightarrow\P(\G)$, and any $Z\in\P(\G)$ the following
equation holds: \hfill

\begin{equation}\label{quad}
\begin{split}
  \theta(Z+W)\times [&\theta(A+U+V+Z)\,\theta(Z-U)\,\theta(Z-V)\\
&-c_1^2c_3^2\,\theta(A+U-V+Z)\,\theta(Z-U)\,\theta(Z+V)\\
&-c_2^2c_3^2\,\theta(A-U+V+Z)\,\theta(Z+U)\,\theta(Z-V)\\
&+c_1^2c_2^2\,\theta(A-U-V+Z)\,\theta(Z+U)\,\theta(Z+V)]=  \\
{}=\theta(A+Z)\times [&\theta(W+U+V+Z)\,\theta(Z-U)\,\theta(Z-V)\\
&-c_1^2c_3^2\,\theta(W+U-V+Z)\,\theta(Z-U)\,\theta(Z+V)\\
&-c_2^2c_3^2\,\theta(W-U+V+Z)\,\theta(Z+U)\,\theta(Z-V)\\
&+c_1^2c_2^2\,\theta(W-U-V+Z)\,\theta(Z+U)\,\theta(Z+V)].
\end{split}
\end{equation}
\end{theo}
To the best of the authors' knowledge equation (\ref{quad}) is a new
identity for Prym theta-functions. For $Z$ such that $\theta(W+Z)=0$
it is equivalent to equation (\ref{cm7d-prym}) with the minus sign chosen.
The second equation of the pair (\ref{cm7d-prym}) can be obtained from
(\ref{51}) considered for the odd case, i.e., for $n+m=1\bmod 2$.
Using theta functional formulas, it can be shown using (\ref{51}) that
equation (\ref{quad}) is equivalent to (\ref{laxdd-prym}).

\section{Abelian solutions of the soliton equations}

In \cite{kr-shio,kr-shio1} the authors introduced a notion of {\it abelian solutions}
of soliton equations which provides a unifying framework the elliptic solutions of these equations and
and algebraic-geometrical solutions of rank 1 expressible in terms of Riemann (or Prym) theta-function.
A solution $u(x,y,t)$
of the KP equation is called {\it abelian\/} if it is of the form
\beq\label{ushio}
u=-2\p_x^2\ln \tau(Ux+z,y,t)\,,
\eeq
where $x$, $y$, $t\in\mathbb C$ and $z\in \mathbb C^n$ are independent variables,
$0\ne U\in\mathbb C^n$, and for all $y$, $t$ the function
$\tau(\cdot,y,t)$ is a holomorphic section of a line bundle $\L=\L(y,t)$ on an
abelian variety $X=\bC^n/\Lambda$, i.e., for all $\l\in\Lambda$
it satisfies the monodromy relations
\begin{equation}\label{monshio}
\tau(z+\l,y,t)=e^{a_\l\cdot z+b_\l}\tau(z,y,t),\quad
\hbox{for some $a_\l\in\mathbb C^n$, $b_\l=b_\l(y,t)\in\mathbb C$}\,.
\end{equation}
There are two particular cases in which a complete characterization of the abelian
solutions has been known for years. The first one is the case $n=1$ of
elliptic solutions of the KP equations. The second case in which a complete characterization of abelian solutions
is known is the case of indecomposable principally polarized abelian variety (ppav).
The corresponding $\theta$-function is unique up to normalization, so that
Ansatz (\ref{ushio}) takes the form $u=-2\p_x^2\ln\theta(Ux+Z(y,t)+z)$.
Since the flows commute, $Z(y,t)$ must be linear in $y$ and $t$:
$
u=-2\p_x^2 \ln\theta(Ux+Vy+Wt+z)\,.
$
Besides these two cases of abelian solutions with known characterization,
another may be worth mentioning.
Let $\G$ be a curve, $P\in\G$ a smooth point, and
$\pi\colon\G\to\G_0$ a ramified covering map such that the curve $\G_0$
has arithmetic genus $g_0 > 0$ and $P$ is a branch point of the covering.
Let $J(\G)=Pic^0(\G)$ be the (generalized) Jacobian of $\G$, let
$Nm\colon J(\G)\to J(\G_0)$ be the reduced norm map as in
\cite{mumford_prym}, and let
$$
X=\ker(Nm)^0\subset J(\G)
$$
be the identity component of the kernel of Nm. Suppose $X$ is compact.
By assumption we have
$$
	\dim J(\G) - \dim X = \dim J(\G_0) = g_0 > 0,
$$
so that $X$ is a proper subvariety of $J(\G)$, and the polarization on $X$
induced by that on $J(\G)$ is {\it not\/} principal.
and define the KP flows on $\overline{Pic^{g-1}}(\G)$ using
the data $(\G,P,\z)$.

In general, since for any $r_0\in\bN$ the space
$\sum_{r\le r_0}\bC\p/\p t_r$ is independent of the choice of $\z$,
for any $\z\in\fm_P\setminus\fm_P^2$ and $0<r<m$
(so in particular for $r=1$), the $r$-th KP orbit of $\F$
is contained in $\F\otimes X$, and so it gives an abelian solution.
Let us call this the {\it Prym-like\/} case.
An important subcase of it is the quasiperiodic solutions of
Novikov-Veselov (NV) or BKP hierarchies.

In the Prym-like case, just as in the NV/BKP case we can put
singularities to $\G$ and $\G_0$ in such a way that $X$ remains
compact, so it is more general than the KP quasiperiodic solutions.
Recall that NV or BKP quasiperiodic solutions can be obtained from
Prym varieties $\Prym(\G,\iota)$ of curves $\G$ with involution $\iota$
having two fixed points.
The Riemann theta function of $J(\G)$ restricted to a suitable translate
of $\Prym(\G,\iota)$
becomes the square of another holomorphic function, which defines the
principal polarization on $\Prym(\G,\iota)$.
The Prym theta function becomes NV or BKP tau function, whose square
is a special KP tau function with all {\it even\/} times set to zero,
so any KP~time-translate of it
\begin{itemize}\itemsep0pt
\item gives an abelian solution of the KP hierarchy
with $n=\dim X$ being one-half the genus $g(\G)$ of $\G$, and
\item defines twice the principal polarization on $X$.
\end{itemize}
A natural question is whether these conditions characterize
the (time-translates of) NV or BKP quasiperiodic solutions.

Hurwitz' formula tells us that in the Prym-like case
$n=\dim(X)\ge g(\G)/2$, where the equality holds only in
the NV/BKP case.  At the moment we have no examples of
abelian solutions with $1<n<g(\G)/2$.

For simplicity we present here a solution to the classification problem of abelian solutions
of the KP equation obtained in
\cite{kr-shio} under an additional assumption on the density of the
orbit $\bC U\bmod\Lambda$ in $X$.
\begin{theo} Let $u(x,y,t)$ be an abelian solution of the KP such that the group $\bC U\bmod\Lambda$ is dense in $X$. Then there exists a unique algebraic curve $\G$ with smooth marked point
$P\in\G$, holomorphic imbedding
$j_0\colon X\to J(\G)$ and a torsion-free rank 1 sheaf $\F\in\overline{{\rm Pic}^{g-1}}(\G)$
where $g=g(\G)$ is the arithmetic genus of $\G$, such that setting with the notation
$j(z)=j_0(z)\otimes\F$
\beq\label{is1}
\tau(Ux+z,y,t)=\rho(z,y,t)\,\widehat\tau(x,y,t,0,\ldots\mid\G,P,j(z))
\eeq
where $\widehat\tau(t_1,t_2,t_3,\ldots \mid \G,P,\F)$ is the KP $\tau$-function
corresponding to the data $(\G,P,\F)$, and  $\rho(z,y,t)\not\equiv0$
 satisfies the condition $\p_U\rho=0$.
\end{theo}
 Note that if $\G$ is smooth then:
\beq\label{is2}
\widehat\tau(x,t_2,t_3,\dots\mid\G,P,j(z))=
\theta\Bigl(Ux+\sum V_it_i+j(z)\Bigm|B(\G)\Bigr)\,
e^{Q(x,t_2,t_3,\ldots)}\,,
\eeq
where $V_i\in\mathbb C^n$, $Q$ is a quadratic form, and $B(\G)$ is the
period matrix of $\G$.
A linearization on $J(\G)$ of the nonlinear $(y,t)$-dynamics for  $\tau(z,y,t)$ indicates the possibility of the existence of integrable systems on spaces of theta-functions of higher level. A CM system is an example of such a system for $n=1$.

Without the density assumption there are examples in which
the KP hierarchy has basically no control
beyond the closure of the orbit, showing the importance of the principal
polarization in a Novikov-like conjecture in which a minimal number of
equation is used to study the nature of $X$. Having this in mind, we
may regard principally polarized Prym-Tjurin varieties \cite{kanev}
as a way to study analogues of Novikov's conjecture.

\end{document}